\newcommand{\arxiv}{1}

\let\modus 2
 \let\modus 1
\newcommand{\ifarxiv}[2]{
\if \modus\arxiv
      #1
\else
      #2
\fi
}

\ifarxiv{
    \documentclass[11pt,a4paper]{amsart}
    \usepackage{amssymb,amsfonts,amsthm}
    \usepackage{math403}

\newcommand{\dedication}[1]{\dedicatory{#1}}
\newcommand{\institute}[1]{\address{#1}}
\newenvironment{pf*}[1]{\begin{proof}[#1]}{\end{proof}}
\newenvironment{pf}{\begin{proof}}{\end{proof}}

}{
    \documentclass[numreferences]{kluwer}
    \usepackage{amsmath-BooLes-kluwer}
}

\usepackage{graphicx}
\usepackage{xcolor}

\newcommand{\commentary}{\relax}
\newcommand{\nocomments}{\renewcommand{\commentary}[1]{\relax}}
\newcommand{\details}{\relax}
\newcommand{\nodetails}{\renewcommand{\details}[1]{\relax}}
\nocomments
\nodetails

\ifarxiv{}{% no this is for the kluwer version
\let\templabelenumi=\labelenumi
\renewcommand{\labelenumi}{\textup{\templabelenumi}}
\let\templabelenumii=\labelenumii
\renewcommand{\labelenumii}{\textup{\templabelenumii}}

\newcommand{\enum}[1]{\textup{(#1)}}

\newcommand{\plref}[1]{{\normalfont \ref{#1}}}

\newcommand{\comment}[1]{\relax}
\newtheorem{theorem}{Theorem}[section]
\newtheorem{lemma}[theorem]{Lemma}
\newtheorem{prop}[theorem]{proposition}

\newtheorem{conjecture}[theorem]{Conjecture}
\newtheorem{problem}[theorem]{Problem}
\newdisplay{dfn}[theorem]{Definition}
\newdisplay{remark}[theorem]{Remark}
\newdisplay{example}[theorem]{Example}
} %end of ifarxiv
\providecommand{\mathscr}{\mathcal} % a priori mathscr is mathcal
\IfFileExists{mathrsfs.sty}{
\usepackage{mathrsfs}}         %  \mathscr produces rsfs fonts
   {%  else try euler fonts
     \IfFileExists{eucal.sty}{
        \usepackage[mathscr]{eucal}} % \mathscr produces Euscript fonts,
   {% else both unavailable, do nothing, mathscr remains mathcal
   }
}
\ifarxiv{}{
\newcommand\ga{\alpha}

\newcommand\gG{\Gamma}
 \newcommand\gD{\Delta}

\newcommand\pl{\partial}%<---abweichende Terminologie

\newcommand\gl{\lambda}
\newcommand\gL{\Lambda}
\newcommand\go{\omega}
\newcommand\gO{\Omega}

\newcommand\gs{\sigma}
\newcommand\gS{\Sigma}

\newcommand{\R}{\mathbb{R}}

\newcommand{\Z}{\mathbb{Z}}

\newcommand\cB{\mathscr{B}}

\newcommand\cD{\mathscr{D}}
\newcommand\cE{\mathscr{E}}

\newcommand\dist{\operatorname{dist}}

\newcommand\Id{\operatorname{Id}}
\newcommand\id{\operatorname{id}}
\newcommand\im{\operatorname{im}}

\newcommand\ind{\operatorname{ind}}

\newcommand\Op{\operatorname{Op}}

\renewcommand\Re{\operatorname{Re}}

\newcommand\sign{\operatorname{sign}}
\newcommand\spec{\operatorname{spec}}

\newcommand\str{\operatorname{str}}
\newcommand\supp{\operatorname{supp}}

\newcommand{\bigsetdef}[2]{\bigl\{ #1 \,\bigm|\, #2\bigr\}}

\newcommand{\scalar}[2]{\langle #1,#2\rangle}

} %end of ifarxiv

\newcommand\Calderon{Cal\-der{\'o}n}
\newcommand\UCP{\textup{UCP}}
\newcommand{\wt}{\widetilde}
\newcommand{\sa}{\textup{sa}}
\newcommand{\half}{{1/2}}
\newcommand{\comp}{\operatorname{comp}}
\newcommand{\dstr}{d_{\textrm{str}}}
\newcommand{\Ell}{\operatorname{Ell}}

\newcommand\Ci{C^\infty}
\newcommand{\CL}{\operatorname{CL}}

\ifarxiv{
\newcommand{\dom}{\operatorname{dom}}{
\renewcommand{\dom}{\operatorname{dom}}}
} %ifarxiv

\newcommand{\dd}[1]{\frac{d}{d#1}}

\newcommand\lla{\langle}
\newcommand\noi{\noindent}

\newcommand\ol{\overline}

\newcommand\rra{\rangle}
\newcommand\tand{\mbox{\ \rm  and }}

\newcommand\too{\longrightarrow}

\newcommand\ii{^{-1}}

\newcommand\<{\subset}
\begin{document}
\ifarxiv{}{
\begin{frontmatter}
\end{frontmatter}

\begin{article}
\begin{opening}
}

\title{The Invertible Double of Elliptic Operators}

\dedication{This short review and work program is dedicated to the memory of Krzysztof
P. Wojciechowski (1953-2008), who was a leader of the investigation of
spectral invariants of Dirac type operators for almost 30 years.}

\ifarxiv{
\author{Bernhelm Booss-Bavnbek}
\address{Department of Science, Systems, and Models\\ Roskilde
University, DK-4000 Ros\-kilde, Denmark} \email{booss@ruc.dk}
\urladdr{http://imfufa.ruc.dk/$\sim$booss}

\author{Matthias Lesch}
\address{Mathematisches Institut,
Universit\"at Bonn, Beringstr. 6, D-53115 Bonn, Germany}
\email{lesch@math.uni-bonn.de, ml@matthiaslesch.de}
\urladdr{http://www.matthiaslesch.de}
\curraddr{Department of Mathematics, University of Colorado at Boulder,
Boulder, CO 80309-0395, USA}

}{
\author{Bernhelm \surname{Boo{\ss}--Bavnbek}\email{booss@ruc.dk}}
\institute{Department of Science, Systems, and Models, Roskilde
University, DK-4000 Ros\-kilde, Denmark}
\author{Matthias \surname{Lesch}\email{ml@matthiaslesch.de, lesch@math.uni-bonn.de}}
\institute{Mathematisches Institut,
Universit\"at Bonn, Beringstr. 6, D-53115 Bonn, Germany\\
current address: Department of Mathematics, University of Colorado at Boulder,
Boulder, CO 80309-0395, USA 
}}


\ifarxiv{
}{
\classification{2000 Mathematics Subject Classification}{
Primary 58J32; Secondary 35J67, 58J50, 57Q20}
}  %ifarxiv

\begin{abstract}
First, we review the {\it Dirac operator folklore} about basic
analytic and geometrical properties of operators of Dirac type on
compact manifolds with smooth boundary and on closed partitioned
manifolds and show how these properties depend on the construction
of a canonical invertible double and are related to the concept of
the {\it \Calderon\ projection}.  Then we summarize a recent
construction of a canonical invertible double for {\it general first order
elliptic differential operators} over smooth compact manifolds with
boundary. We derive a natural formula for the \Calderon\ projection
which yields a generalization of the famous Cobordism Theorem. We
provide a list of assumptions  to obtain a continuous
variation of the \Calderon\ projection under smooth variation of the
coefficients. That yields various new spectral flow theorems.
Finally, we sketch a research program for confining, respectively
closing, the last remaining gaps between the geometric Dirac
operator type situation and the general linear elliptic case.
\end{abstract}

\ifarxiv{\maketitle}{
\end{opening}}
\setcounter{section}{-1}
\section{Introduction}

This paper reviews our recent results, obtained jointly with Chaofeng Zhu \cite{BooLesZhu:CPD},
about basic analytical properties of elliptic operators on compact
manifolds with smooth boundary. Furthermore, we outline a research
program for confining, respectively
closing, the last remaining gaps between the geometric Dirac
operator type situation and the general linear elliptic case.

Our main results are 
\begin{itemize}
\item to develop the
basic elliptic analysis in full generality, and not only for the generic case of 
operators of Dirac type in product metrics (i.e., we assume neither constant coefficients
in normal direction nor symmetry of the tangential operator); 
\item to give an analytical proof of the
cobordism invariance of the index in greatest generality; and 
\item to prove the continuity of the \Calderon\ projection and 
of related families of global elliptic boundary
value problems under parameter variation. 
\end{itemize}

Most analysis of geometrical and physical problems involving a Dirac operator $A$ on a compact manifold $M$ with smooth 
boundary $\partial M$ acting on sections of a (complex) bundle $E$ seems to rely on quite a few basic facts which are part 
of the shared folklore of people working in this field of global analysis. See, 
e.g., \cite{BooWoj:EBP} for 
\begin{itemize}
\item the weak inner unique continuation property (also called weak UCP to the boundary), i.e., there are no nontrivial elements in the null space $\ker(A)$ vanishing at the boundary of $M$; 

\item the existence of a suitable elliptic invertible continuation $\wt A$ of $A$, acting on sections of a vector bundle over the closed double or another suitable closed manifold $\wt M$ which contains $M$ as submanifold; this yields a Poisson type operator $K_+$ which maps sections over the boundary into sections in $\ker A$ over $M$; and a precise \Calderon\ projection $C_+$\,;

\item the existence of a self-adjoint regular Fredholm extension of any total (formally self-adjoint) Dirac operator $A$ in the underlying $L^2$-space with domain given by a pseudodifferential boundary condition; that actually is equivalent to the Cobordism Theorem asserting a canonical splitting of the tangential operator $B = B^+ \oplus B^-$ with $\ind(B^+) = 0$; 
\end{itemize}

\noi and Nicolaescu \cite[Appendix]{Nic:GSG} and Boo{\ss}, Lesch and Phillips \cite{BooLesPhi:UFOSF} for 
\begin{itemize}
\item the continuous dependence of a family of operators, their associated \Calderon\ projections, and of any family of well-posed (elliptic) boundary value problems on continuous or smooth variation of the coefficients.
\end{itemize}

We will show that these results, with slight modifications, do hold for arbitrary first order
elliptic operators.

We are fully aware that Dirac operators are attractive because they have an
obvious geometrical meaning like the Laplace operator, probably the most important elliptic
operator, for physicists as well as mathematicians. In the present paper we show that the above
mentioned results for the Dirac operator do hold, with slight modifications, for arbitrary first order
elliptic operators. As a matter of fact, the methods applied can be extended to general elliptic
systems (of higher order) in a very natural way. However, in this paper we shall restrict ourselves
to first order systems, mostly for the ease of presentation and the clarity of the essential
constructions and proofs.

It is therefore not unreasonable to ask 
\textit{How special are operators of Dirac type compared to 
arbitrary linear first order elliptic differential operators?}
Our treatment of that question here may be used as a guide for addressing corresponding questions
for general elliptic systems of higher order. Roughly speaking, our message is that the basic
geometric aspects of geometrically defined operators like the Dirac operator (and the Laplace
operator, when extending to higher order) have their root in the very ellipticity, i.e., in the
symmetry properties of the principal symbol. The beautiful underlying metric properties leading to
the definition of these operators make proofs of the basic geometric facts easier but are, as we
shall show, dispensable.

We were led to our deformation question by a variety of mathematical and physical motivations. 
We mention 
\begin{itemize}
\item the continuing interest in contact manifolds and
CR-structures and the corresponding tangential CR-complexes (see, e.g.,
the classic monograph Boggess \cite{Bog:CRM} or the recent Ponge \cite{Pon:NRI});
\item stability questions in Electrical Impedance Tomography (see, e.g., the original
problem in \Calderon\ \cite{Cal:IBV} and the recent Kenig and Sj{\"o}strand
\cite{KenSjoUhl:CPP}, dealing, though, with
the robustness of the Dirichlet-to-Neumann operator for elliptic second order equations instead
of the \Calderon\ projection for first order operators dealt with in this Note);
\item deformations in Quantum Gravity which tear and/or burst classical field theories
(see the visions disseminated, e.g., in Nielsen and Ninomiya \cite{NieNin:LBS});
\item but most of all our own curiosity about the extent to which simple geometric properties
not only predetermine and guide analytic investigation but also pin down the results.
\end{itemize}
 
It has been known for half a century that, e.g., the K-groups of
spin manifolds are generated by the index classes of Dirac operators
up to torsion. For concrete calculations, however, many technical
arguments depend on constructions which work only for geometrically
defined operators of Dirac type and under the additional assumption
that all metric structures are product near the boundary $\partial M$ of
the underlying smooth compact Riemannian manifold $M$.

One of such technical devices is the invertible double, proved
in Boo{\ss} and Wojciechowski \cite{BooWoj:EBP}.
For the precise formulation and a sketch of proof see below 
Proposition \ref{p:inv-dirac}.

For more general elliptic operators, as arising from first-order
deformations, an invertible
elliptic extension is often assumed ``for convenience''. Actually,
for Dirac type operators in the non--product case, one can still
extend the collar a bit and deform to the product situation. The
resulting operator will still be invertible; however it will neither
be an exact double, nor will it be canonical.

Although a geometric invertible double is not available in general,
in this Note we shall show that there is always a nice boundary
value problem which provides an ``invertible'' double and that 
important properties, previously established rigorously only for
Dirac type operators remain valid for general elliptic operators.

%As always with so-called {\it technicalities}, it requires some effort 
%to distinguish clearly between what is known, what is thought to be known but may not be
%so simple to achieve, what is new, and what remains to be done. That dictates the 
%structure of our research and review article.

The paper ist organized as follows:

In Section \ref{s:folklore}, we summarize the Dirac operator {\em folklore} with emphasis on
the product property, weak inner UCP, the precise invertible double, and 
a sketch of the geometric role of the \Calderon\ projection. 
Some of the results are counter-intuitive in spite of their basic and fundamental character.
E.g., the local solvability of elliptic equations is well-known from classical theory,
but it remains a surprise that, e.g., the precise double of the 
Cauchy--Riemann operator (obtained by twisting the complex line bundle, see below) 
is well defined, 
has smooth coefficients and is invertible - without subtracting projections on 
original or arising kernels.

In Section \ref{s:invertible}, we present the first main result of this article, namely
the construction of a precise invertible double for any first order elliptic differential operator, 
satisfying weak inner UCP, respectively, an invertible double after subtracting the
projection onto the inner ({\em ghost}) solutions. The novelty of our approach lies in
the canonical character of the construction - in difference to the ingenious ideas of 
Seeley \cite{See:SIB}, \cite{See:TPO} of the late 1960s which also provided an 
invertible double,
but involved extensions and choices which excluded to follow, e.g., the parameter dependence
and neither yielded the Lagrangian property of the Cauchy data spaces in the case of 
symmetric coefficients.

In Section \ref{s:appl}, we present the other main results of this article, namely various applications
of our construction of the invertible double. Surprisingly, it turns out that
the investigation of the
mapping properties of the induced Poisson operators and \Calderon\ projections is by no means
straightforward. It may be worth recalling the decisive role of the socalled
{\em Atiyah--Patodi--Singer} spectral projections of the induced symmetric
tangential operator over the boundary in the Dirac operator case. In our general case, 
that nice tool must be replaced by {\em sectorial projections}.
Nevertheless, we can in this Section \ref{s:appl} establish
\begin{itemize}
\item the Lagrangian property of the Cauchy data space for formally self-adjoint coefficients;
\item a Cobordism Theorem for any elliptic operator bounding a first order elliptic differential operator; and
\item a couple of theorems analyzing the dependence of the \Calderon\ projection on the input data.
\end{itemize}

We emphasize that all these results are stated for {\em general} first order elliptic differential operators. 
It is neither assumed that the operator is of Dirac type nor is product structure near the boundary assumed.

The proofs are intricate. We only give sketches of the proofs and made full-length proofs available at arXiv, Boo{\ss}, Lesch and Zhu
\cite{BooLesZhu:CPD}.

%To emphasize the envisaged general interest and geometrical meaning of our new results 
%(achieved jointly with Chaofeng Zhu of the Chern Institute of Nankai University, China), 
%we embedded the results in the present review and work program note.

Three years ago, before our results were achieved, the first author made a ``poll" at a 
conference of experts in global analysis about the correctness of our at that time only conjectures: 
about one half of the people present at that meeting thought the claims were more or less clear and almost 
proved already in the late 60s or early 70s. 
The other half %(perhaps the half with a better understanding of the power of geometric reasoning associated 
%to the concept of Dirac operators and product structures and the risks when loosing that power by deformation) 
doubted the claims and would bet on counter-examples.

In Section \ref{s:gaps}, we explain why we do not consider the reached results for optimal;
what difficulties must be overcome; and what ideas might turn out to be worth following.
In particular, it seems to us that a much better understanding of the analysis and geometry of 
sectorial projections is mandatory for further work on the mapping properties; and that 
the time perhaps has come for new approaches towards UCP, namely by focusing on 
weak inner UCP.% - after the 70 years
%of post-Carleman stagnation in that field.

\section{Dirac operator {\em folklore}}\label{s:folklore}

%%\subsection{Product property of operators of Dirac type}

%%\noi 1.1 \textit{Product property of operators of Dirac type}.
\paragraph*{\textup{1.1} \textit{Product property of operators of Dirac type.}}
Let $M$ be a smooth compact oriented Riemannian manifold of dimension
$m$ (with or without boundary) and let $E\to M$ be a Hermitian
vector bundle of Clifford modules
with the Clifford multiplication ${\mathbf c}:\Gamma^\infty(M;TM\otimes E)\to \Gamma^\infty(M;E)$.
We recall that any choice of a smooth connection (covariant derivative)
$\nabla:\Gamma^\infty(M;E)\to\Gamma^\infty(M;T^*M\otimes E)$ defines a
{\em (total) operator of Dirac type}
$ A:={\mathbf c}\circ \nabla:\Gamma^\infty(M;E)\too\Gamma^\infty(M;E)$,
acting on the space $\Gamma^\infty(M;E)$ of smooth sections
under the Riemannian identification of the bundles $TM$ and $T^{\ast}M$.

In local coordinates we have
\begin{equation}
A=\sum_{i,j=1}^{m} g^{ij}\mathbf{c}(\frac{\pl}{\pl x_i})\frac{\partial
}{\partial x_{j}}+\text{ zero order terms}. \ifarxiv{\label{e:dirac}}{\label{ie:dirac}}       %
\end{equation}
It follows at once that the principal symbol $\sigma_{1}(D)(p,\xi)$
is given by Clifford multiplication by $i\xi$, so that any operator of Dirac
type is elliptic with symmetric principal symbol. 
Denoting by $A^t$ the formal adjoint of $A$ we have 
{\em Green's formula}
      \begin{equation}
     ({A} s,s') - (s,{A^t} s') = -\int_{\pl M}
     \langle J (s|_{\pl M}),s'|_{\pl M}\rangle ,\quad
s,s'\in \Gamma^\infty(M;E).
\ifarxiv{\label{e:green}}{\label{ie:green}}
     \end{equation}
Here $J:={\mathbf c}({\bf n}): E|_{\pl M} \to E|_{\pl M}$ denotes the unitary
bundle isomorphism given by Clifford multiplication by the
inward unit tangent vector with $J^2=-I$.

If the connection $\nabla$ is \emph{compatible} with Clifford multiplication 
(i.e. $\nabla\mathbf{c}=0$) and \emph{unitary} (i.e. Leibniz' rule
$X\langle s, s'\rangle=\langle \nabla_X s,s'\rangle+\langle  s,\nabla_X s'\rangle$ holds
for $s,s'\in\Gamma^\infty(M;E), X\in\Gamma^\infty(TM)$), 
then the operator $A$ itself becomes formally self--adjoint.

\details{It is indeed true that Green's formula does not depend on any specific
choice (compatible, unitary) of the connection. Integration by parts affects only
the first order part of the connection}

Let $\gamma_5$ denote the global section of $\operatorname{Hom}(E,E)$
defined locally by $\gamma_5:= {\mathbf c}(e_1)\dots {\mathbf c}(e_{m})$
(for a positively oriented orthonormal local frame).
If $m$ is even, e.g., $m=4$ as in many physics applications,
$E$ splits into subbundles $E^{\pm}$\,. They are spanned by the eigensections of $\gamma_5$
corresponding to the eigenvalue $\pm 1$, if $m$ is divisible by
4, or $\pm i$ otherwise. The Clifford multiplication $J$
switches between $E^\pm|_{\pl M}$ and $E^\mp|_{\pl M}$. If $\nabla$ is compatible
and unitary\footnote{This condition can certainly be somewhat relaxed, e.g. for the splitting
of the Dirac operator one just needs that the decomposition $E=E^+\oplus E^-$ is parallel with
respect to $\nabla$.} the Dirac operator splits correspondingly into components
     $
     {A}=\bigl(\begin{smallmatrix}
     0 & A^-\\ {A}^+ & 0
     \end{smallmatrix}\bigr)
     $
such that the {\em right chiral (half) Dirac operator}
${A}^+: \Gamma^\infty(M;E^+)\to\Gamma^\infty(M;E^-)$ is formally adjoint to
${A}^-:\Gamma^\infty(M;E^-)\to\Gamma^\infty(M;E^+)$.

From \eqref{e:dirac} we derive a product property which
distinguishes operators of Dirac type from {\em general} elliptic
differential operators of first order.

%% \medskip

\begin{lemma} \label{l:product}
Let $\Sigma$ be a closed hypersurface of $M$ with orientable
normal bundle. Let $x$ denote a normal variable with fixed orientation such
that a bicollar neighborhood ${N}$ of $\Sigma$ is parametrized by
$[-\varepsilon,+\varepsilon]\times\Sigma$. Then any operator of Dirac type can
be rewritten in the form
\begin{equation}
{A}|_{{N}}=\mathbf{c}(\frac{\partial}{\partial x})\left(  \frac{\partial}{\partial x}%
+B_{x}+C_{x}\right)  , \ifarxiv{\label{e-product}}{\label{ie-product}}%
\end{equation}
where $B_{x}$ is a self-adjoint elliptic operator on the parallel hypersurface
$\Sigma_{x}$, and $C_{x}:E|_{\Sigma_{x}}\rightarrow E|_{\Sigma_{x}}$ is a
skew-adjoint operator of $0$th order, i.e., a skew-symmetric bundle homomorphism.\bigskip
\end{lemma}
We shall call the operator $B_0+C_0$ the {\em tangential} component of $A$ in direction $\Sigma$.
Its principal symbol is symmetric.
%% \medskip
\details{To avoid further confusion of ML:

in this lemma we neither assume $A$ to be symmetric nor do we assume anything about
the connection. Note that choosing a compatible and unitary connection affects only
the $0$th order terms of the operator. So since for a compatible and unitary connection 
$A$ is self--adjoint, indeed $A-A^t$ is $0$th order! By the same reasoning the tangential
operator has self--adjoint leading symbol. This explains the decomposition of the tangential
operator into a self--adjoint first order operator and a skew--adjoint $0$th order operator.
I think ML has now understood it, too.}

If the Riemannian
metric of $M$, the Hermitian structure of $E$ and the connection $\nabla$ are product near the
boundary, it follows from Lemma \ref{l:product} (re-writing $B_0+C_0=:B$)
that each operator of Dirac type takes the form
\begin{equation}\ifarxiv{\label{e:bord-product}}{\label{ie:bord-product}}
A=J(\dd x + B)\quad\text{close to $\partial M$,}
\end{equation}
with $J,B$ independent of the
normal variable $x$ and $B$ a first order
elliptic differential operator on $\partial M$. $B$ is not necessarily self--adjoint
but has self--adjoint leading symbol. 
If the connection $\nabla$ is compatible and unitary then $A$ and $B$ are (formally) self--adjoint
and this then implies $JB=-BJ$.
\details{Indeed the formula \eqref{e:bord-product} has
nothing to do with symmetry of $A$ resp. compatibility of $\nabla$. 
If all structures including the connection are product, then near a boundary point we may choose
coordinates $x=x_0, x_1,...,x_n$ ($x=x_0$ denotes the normal coordinate)
such that $g_{0j}=\delta_{0j}$ and such that the connection form of $\nabla$
takes the form $\sum_{j\ge 1} A_j(x_1,...,x_n) dx_j$.
The latter is just the definition of $\nabla$ being in product form. Of course this
additional assumption on $\nabla$ does not follow from the product assumptions on the metrics.
Note however, that for the spin connection on a spin manifold it \emph{follows} from the
product assumptions on the metrics that the spin connection is in product form, too.

Anyway it then follows that
\[\begin{split}
       D&=\mathbf{c}(\partial_0)\frac{\partial}{\partial x_0}+\sum_{i,j\ge 1} g^{ij} \mathbf{c}(\partial_i)\nabla_{\partial_j}\\
          &= \mathbf{c}(\partial_0)\bigl(\frac{\partial}{\partial x_0}+
         \sum_{i,j\ge 1} g^{ij} \mathbf{c}(\partial_0)\mathbf{c}(\partial_i)\nabla_{\partial_j}\bigr)\\
          &= J\bigl(\frac{\partial}{\partial x}+B\bigr)
\end{split}
\]
with $B=\sum_{i,j\ge 1} g^{ij} \mathbf{c}(\partial_0)\mathbf{c}(\partial_i)\nabla_{\partial_j}$
and $J=\mathbf{c}(\partial_0)$. The Clifford relations immediately imply that for 
$JB+BJ=0$ to hold it is sufficient that $\mathbf{c}(\partial_0)$ is parallel. This does not
follow from the product structure but does probably not need the full strength of 
a compatible and unitary connection.
}

%% \subsection{Unique Continuation Property (\UCP) for operators of Dirac type}
%%\noi 1.2 \textit{Unique Continuation Property (\UCP) for operators of Dirac type}.
\paragraph*{\textup{1.2} \textit{Unique Continuation Property (\UCP) for operators of Dirac type.}}
%% \begin{e-definition}\ifarxiv{\label{d:ucp}
Let $A$ be a linear or non--linear operator, acting on functions or sections of a bundle $E$ over a
compact or non--compact connected manifold $M$.

\noi (a) The operator $A$ has the {\em strong} \UCP\ if any solution $u$ of the equation $Au=0$ has the following property: if $u$ vanishes at a point $p\in M$ with all its derivatives, then it vanishes on the whole of $M$.

\noi (b) The operator $A$ has the {\em weak} \UCP\ if any solution $u$ of the equation $Au=0$ has the following property: if $u$ vanishes on a nonempty open subset $\gO$ of $M$, then it vanishes on the whole of $M$.

\noi (c) Let $M$ be a compact manifold with boundary $\pl M$.
Then the operator $A$ has the {\em weak inner} \UCP\ (also called {\em weak \UCP\ ``across the boundary"})
if any solution $u$ of the equation $Au=0$ has the following property:
if $u$ vanishes on $\pl M$, then it vanishes on the whole of $M$.

%% \end{e-definition}

Let $A$ be an arbitrary elliptic differential
operator $A$ of {\em first} order on a closed partitioned manifold $M=M_-\cup_{\Sigma} M_+$\,, where $\Sigma$ is a hypersurface.
Elliptic regularity and Green's Formula imply (see, e.g., \cite[Lemma 12.3]{BooWoj:EBP}):

\begin{lemma}\label{l:partitioned}
Any $u_+\in\Gamma^\infty(M_+; E|_{M_+})$ with $Au_+=0$ and $u_+|_{\Sigma}=0$ can be continued to a smooth solution $u$ for the operator $A$ over the whole
manifold $M$ by setting $u :=(u_+,0)$.
\end{lemma}

The preceding Lemma implies that weak \UCP\  can be reformulated for linear elliptic differential operators of first order as \UCP\  across {\em any} hypersurface. More precisely, the operator $A$ has the weak \UCP\ if any solution $u$ of the equation $Au=0$ has the following property: if $u$ vanishes on a hypersurface $\Sigma$, then it vanishes on the whole of $M$.
In particular, weak \UCP\  implies weak {\em inner} \UCP\  for linear elliptic differential operators of first order.

Finally, by the same argument we see that weak \UCP\  implies \UCP\  across any single connected component of the
boundary of $M$. More precisely, let $u$ be any solution of the equation $Au=0$ with $u|_{\Sigma_1}=0$ where $\Sigma_1$ is one connected component of the boundary $\Sigma$ of a compact connected manifold $M$. Then by weak \UCP\  it vanishes on all other components of $\Sigma$ and, in particular, it vanishes on the whole of $M$. Like so many other features of complex analysis, this tunneling property associated to weak \UCP\  is a bit counter--intuitive.

Weak \UCP\  is one of the basic properties of any operator of Dirac type. This can be seen by applying a hard result, obtained in 1956 independently by N. Aronszajn and H. O. Cordes for
linear scalar elliptic operators of second order with smooth coefficients and with real principal symbol, to the Dirac Laplacian.

While the Aronszajn--Cordes Theorem yields strong \UCP, there is also a direct proof on the level of the Dirac operator,
exploiting Lemma \ref{l:product}, but yielding only weak \UCP.

\begin{theorem}\label{t:ucp-symmetric}
Let $M$ be a smooth compact connected Riemannian manifold and  $A$ an elliptic differential operator
of first order acting on sections in a Hermitian bundle. We assume that
the tangential operator has symmetric principal symbol on every hypersurface $\Sigma$. Then $A$
satisfies weak \UCP.
\end{theorem}

%\ifarxiv{
%\begin{proof}[Sketch of proof]}{
\begin{pf*}{Sketch of proof}
%}
Let $u$ be a solution of $Au=0$ which vanishes on an open nonempty set $\Omega$. We assume $\Omega=M\setminus \supp u$, i.e., maximal, namely the union of all open subsets on which $u$ vanishes.

First we localize and convexify the situation: since $M$ is connected, to prove that $\Omega=M$ it suffices to show that $\Omega$ is closed. So assume that $\partial \Omega=\ol\Omega\setminus \Omega$ is nonempty. Choose a $p\in \Omega$ whose
distance from $\partial \Omega$ is less than the injectivity radius of $M$.
Then we choose $p_0\in\partial\Omega$ with $\dist(p,p_0)=\dist(p,\partial\Omega)$.
In other words the open ball around $p$ with radius $r:=\dist(p_0,p)$ is contained in $\Omega$, but $p_0\in\partial\Omega\subset\supp u$.

This construction provides us with a family of concentric hyperspheres $S_{p,x}$ of radius $x+r$. We
fix an angular region by choosing $T> 0$ with $T$ sufficiently small and inner radius $r$, ranging from the hypersphere $S_{p,0}$ which is contained in $\ol\Omega$, to the hypersphere $S_{p,T}$ which cuts deeply into $\supp u$, if $\supp u$ is not empty.

To conclude the localization, we replace the solution $u|_{ \{S_{p,x}\}_{x\in [0,T]}}$ by a cutoff $v(x,y):=\phi(x) u(x,y)$
with a smooth bump function $\phi$.

Now we establish a Carleman inequality: for all $T >0$ sufficiently small and all $R>0$ sufficiently large we have
\begin{multline}\ifarxiv{\label{e:carleman}}{\label{ie:carleman}}
    \int_0^T \!\int_{S_{p,x}}\!\! e^{R(T-x)^2} \| v(x,y)\|^2 dy\,dx\\
     \le  \frac{2}{R}\int_0^T \!\int_{S_{p,x}}\!\! e^{R(T-x)^2} \| A v(x,y)\|^2 dy\,dx,
\end{multline}
for all arbitrary smooth sections $v$ (not necessarily a solution) with $\supp v\subset \{S_{p,x}\}_{x\in [0,T]}$\,.

We apply \eqref{e:carleman} to our cutoff section $v(x,y)$ which by construction coincides with the solution $u$ for $x\le 4/5T$ and vanishes identically for $x\ge 9/10T$. Elementary integral inequalities yield, that the solution $u$ must vanish in the whole annular region $0\le x\le T/2$\,. This contradicts $p_0\in\supp u$, which proves the theorem.
%\ifarxiv{
%\end{proof}}{
\end{pf*}
%}

Details of the proof are given in \cite[Chapter 8]{BooWoj:EBP} and further elaborated, e.g., in
Bleecker and Boo{\ss} \cite{BleBoo:SIO} and Boo{\ss}, Marcolli and Wang \cite{BooMarWan:WUCP}.
Already Nirenberg \cite[Sections 6--7, in particular the proof of his inequality (7.11)]{Nir:LLP}
pointed to the decisive role of the symmetry condition for deriving Carleman type inequalities.

%% \subsection{The invertible double for operators of Dirac type}
%%\noi 1.3 \textit{The invertible double for operators of Dirac type}.
\paragraph*{\textup{1.3} \textit{The invertible double for operators of Dirac type.}}
After these preparations, we recall the invertible double construction from
\cite[Chapter 9]{BooWoj:EBP}.

\begin{prop}\label{p:inv-dirac}
Let $M$ be a smooth compact connected Riemannian
manifold with boundary and let
$A$ be a Dirac type operator on $M$ acting between sections of the
Hermitian vector bundle $E$. Assume that all structures are product
near the boundary. Then $A$ and $-A$ can be glued together to obtain
an \emph{invertible} elliptic operator $\tilde A=A\cup_{\partial M}
(-A)$ on the closed double $\tilde M$.
\end{prop}

%\ifarxiv{
%\begin{proof}[Sketch of proof]}{
\begin{pf*}{Sketch of proof}
%}
For
simplicity we shall sketch the proof only for the total Dirac type operator and only in the
formally self-adjoint case.
The non-symmetric case and the chiral component operators
require many, mostly notational modifications, see \cite[Chapter 9]{BooWoj:EBP}.

Our starting point is the product form \eqref{e:bord-product}.
The unitary map sending a section $f$ over $[0,\varepsilon)\times \partial M$
to $Jf(-.)$ over $(-\varepsilon,0]\times \partial M$ conjugates
$A$ and $-A$. Hence if we use $J$ as a clutching function for the bundle
$E$ then indeed $\tilde A=A\cup_{\partial M}(-A)$ is an elliptic operator
with smooth coefficients on the double. Sections in $\dom(\tilde A)$ can
be viewed as pairs $(f_+,f_-)$, $f_\pm=f|_{ M_\pm}$ such that $f_-|_{\pl M}=J(f_+|_{\pl M})$.
Here $M_\pm$ denote the two different copies of $M$ in $\tilde M$.

To prove that $\widetilde A$ is invertible assume that $\tilde A f=0$.
We then have $Af_\pm=0$. Green's formula gives
\begin{equation*}
%%  \begin{split}
     \langle f_{\partial M},f_{\partial M}\rangle
%% = \langle J (f|_{M_-})|_{\partial M},
%% (f|_{M_+})|_{\partial M}\rangle
  =\langle A f_-,f_+\rangle -\langle  f_-,A f_+\rangle=0.
%% \end{split}
\end{equation*}
Thus $f|_{\partial M}=0$. The weak \UCP\  of Dirac operators (Theorem \ref{t:ucp-symmetric})
then implies $f=0$.
%\ifarxiv{
%\end{proof}}{
\end{pf*}
%}

An illustration of the construction is given in \cite[Chapter 26]{BooWoj:EBP} in the simplest
possible two--dimensional case, namely for the Cauchy--Riemann operator on the disc.

%% \subsection{The geometric role of the \Calderon\ projection}

%%\noi 1.4 \textit{The geometric role of the \Calderon\ projection}.
\paragraph*{\textup{1.4} \textit{The geometric role of the \Calderon\ projection.}}
The concept of an invertible double was used already 40 years ago
to show that the
{\em Calder{\'o}n projection} $C_+$, i.e., the projection of $L^2(\partial M; E|_{\pl M})$ onto
the space of Cauchy data $N_+^0:= \{u|_{\partial M}\in
L^2(\partial M; E|_{\pl M})\mid Au=0\}$
is a pseudo--differential operator
(Seeley \cite{See:SIB}, \cite{See:TPO}) for any elliptic differential operator
of first order. Various choices, however, entered into the original construction of the invertible
double, while the preceding construction for Dirac type operators and product metrics close
to the boundary is canonical. That provides
a formula for $C_+$ in terms of
$A$ such that the Lagrangian property $(N_+^0)^\perp
=JN_+^0$ is implied, see \cite[Corollary 12.6]{BooWoj:EBP}, and mapping properties and
dependencies on the data become transparent, see \cite{BooLesPhi:UFOSF}.
Only recently, the authors
of this short review were able
to prove similar results for general elliptic differential operators of first order, in
collaboration with C. Zhu \cite{BooLesZhu:CPD} (see Theorems
\ref{t:calproj}, \ref{t:cobord}, and \ref{t:param} below).

To give the reader the taste of the geometric meaning of $C_+$ we recall
a few results for spectral invariants of Dirac type operators
on manifolds with smooth boundary and on closed partitioned manifolds.

\noi (A) The index of a well-posed boundary value problem $A^+_{P}$ for $A^+$ defined by
a pseudo-differential projection $P$ with the same principal symbol as $C_+$ is given by the
relative index of $P,C_+$, that is the index of the Fredholm operator
$PC_+:\im C_+\to \im P$. On a closed partitioned manifold the index of
a (chiral) Dirac operator $A^+$ can be identified with the index of the
Fredholm pair of Cauchy data spaces along the partitioning hypersurface $\Sigma$.

For various extensions of this {\it Bojarski Formula} to the
subelliptic case and for recent implications for the Atiyah-Weinstein
conjecture see Epstein \cite{Eps:CRI} and references given there.

An interesting feature of the \Calderon\ projection for the Euclidean Dirac operator on 
the 4--ball is that it is $\gamma_5$ invariant, defines a self-adjoint boundary problem, 
and that the corresponding domain is gauge--invariant. As shown in 
Boo{\ss}-Morchio-Strocchi-Wojciechowski \cite{BooMorStrWoj:GCA}, 
this property of the \Calderon\ projection refutes the common claim of Quantum Chromodynamics according to which so-called {\em naturality} (i.e., gauge invariance of the domain 
and self-adjointness and $\gamma_5$ invariance of the boundary condition) implies 
chiral anomaly. The claim was suggested in Ninomiya and Tan \cite{NinTan:AAI}. 
It was based on a (here) misleading property of
the Atiyah-Patodi-Singer boundary condition, namely the non-vanishing index of the chiral problem in general. Actually, this {\em chiral anomaly} would be an obstruction for zeta-function regularization of the determinant, as explained in \cite{BooMorStrWoj:GCA}. Happily, 
the \Calderon\ projection defines a boundary condition which gives vanishing index and, in fact, vanishing kernel and cokernel.

\noi (B) The spectral flow of a curve of (total) Dirac operators with continuously varying
 connections over a closed partitioned manifold equals the Maslov index of the corresponding curves of Cauchy data spaces.

\noi (C) The $\zeta$--determinant of a well-posed self-adjoint boundary value problem $A_P$ equals the Fredholm
determinant of a canonically associated operator over the boundary up to a constant which can be identified with the
$\zeta$--determinant of $A_{C_+}$\,.

For proofs, we refer to \cite{BooWoj:EBP} for (A), to Nicolaescu \cite{Nic:GSG} for (B), and to Scott and Wojciechowski \cite{ScoWoj:ZDQ} for (C).

% main text

\section{Invertible double for general first order elliptic
operators}\label{s:invertible}

Let $M$ be a smooth compact connected Riemannian manifold with
boundary and let $A:\Gamma^\infty(M;E)\longrightarrow
\Gamma^\infty(M;F)$ be a first order elliptic differential operator
acting between sections of the Hermitian vector bundles $E,F$.
As above we separate variables in a collar $U$ of the boundary and
write
\begin{equation}\ifarxiv{\label{eq:operator-collar}}{\label{ieq:operator-collar}}
\begin{split}
   A|_U&= J_0\bigl(\frac{d}{dx}+B_0\bigr)+C_1x+C_0,\\
  A^t|_U&= \bigl(-\frac{d}{dx}+B_0^t\bigr)J_0^*+\wt C_1x+\wt C_0,
\end{split}
\end{equation}
with
bundle morphisms $J_0,C_0,\wt C_0$\,; $B_0$ a first order elliptic
differential operator on $\partial M$; and $C_1,\wt C_1$ first order
differential operators on $U$. Put
%% \begin{equation}
$\tilde A:=A\oplus (-A^t)$,
%% \end{equation}
acting on sections of $E\oplus F$. $A^t$ denotes the formal adjoint
of $A$. We choose a bundle morphism $T\in
\operatorname{Hom}(E|_{\partial M},F|_{\partial M})$ and impose the
boundary condition $({f_+},{f_-})\in\operatorname{dom}(\tilde
A_T):\Leftrightarrow
       {f_-}|_{\partial M}=T {f_+}|_{\partial M}
\Leftrightarrow ({f_+}|_{\pl M}, {f_-}|_{\pl M})\in\ker \begin{pmatrix}
-T &\Id \end{pmatrix}.$

The two most important cases are $T:=(J_0^*)^{-1}$ and, if
$J_0=-J_0^*$, $T:=J_0 |J_0|^{-1}$. In both cases the endomorphism
$J_0^*T$ is positive definite.

\begin{theorem}\label{t:invdoub} Assume that $J_0^*T$ is positive definite.
Then
%\begin{enumerate}

\textup{1.} $\tilde A_T$ is a realization of a local elliptic boundary
condition (in the classical sense of {\v S}apiro-Lopatinski{\v i}), hence
$\tilde A_T$ is a Fredholm operator with compact resolvent.

\textup{2.} $\ker \tilde A_T=Z_0(A)\oplus Z_0(A^t)\ \tand\ \operatorname{coker} \tilde A_T
\cong \ker \tilde A_T^*=Z_0(A^t)\oplus Z_0(A)$.
%\end{enumerate}
\end{theorem}
\begin{remark}
Here $Z_0(A)=\bigl\{ u\in L^2(M,E)\,|\, Au=0, u|_{\partial
M}=0\bigr\}$ denotes the space of ``ghost solutions''.
By elliptic regularity it is easy to see that $Z_0(A)$ is a finite--dimen\-sional
subspace of $\Gamma^\infty(M;E)$ and hence does not depend on the
choice of a Sobolev regularity for $u$. $Z_0(A)=\{0\} $ if and only
if weak inner \UCP\  holds for $A$. While weak \UCP\  can be proved for
Dirac type operators in various ways (see above Section 1.4), 
it is generally believed
that weak inner \UCP\ does not hold for all first order elliptic differential
operators,
and it is open whether weak \UCP\  for $A$ implies weak
\UCP\  for $A^t$, as conjectured by L. Schwartz 
\cite{Sch:EDP}, cf. below Section 4.2 for all that.

%% \footnote{Ecuaciones diferenciales parciales elipticas, deuxi\`eme \'ed.,
%% {Revista colombiana de matematicas} 13, Bogota, 1973.}

If the operator $A$ is formally self--adjoint, then for $T:=J_0
|J_0|^{-1}$ the double $\tilde A_T$ is self--adjoint.
\end{remark}

%\ifarxiv{
%\begin{proof}[Sketch of proof of Theorem \ref{t:invdoub};
%for details cf. \mbox{\cite[pp. 6--12, 20--23]{BooLesZhu:CPD}}]}{
\begin{pf*}{Sketch of proof of \textup{Theorem \ref{t:invdoub}};
for details cf. \mbox{\cite[pp. 6--12, 20--23]{BooLesZhu:CPD}}}
%}
As explained above the boundary condition for $\wt A_T$ is
given by the bundle homomorphism
$P(T) = \begin{pmatrix}  -T & \Id \end{pmatrix}$.

From \eqref{eq:operator-collar} we see that the tangential
operator of $\wt A$ has leading symbol $b_0\oplus -(J_0^t)\ii b_0^* J_0^t$,
$b_0:=\sigma^1_{B(0)}$. Consequently the positive spectral projection of
$b_0\oplus -(J_0^t)\ii b_0^* J_0^t$ is given by $P_+(b_0)\oplus (J_0^t)\ii
P_-(b_0^*) J_0^t$.
For an endomorphism $b$ of a finite--dimensional
vector space $P_\pm(b)$ denotes the spectral projection corresponding
to a closed contour encircling all eigenvalues $\lambda$ with $\Re \lambda\ge 0 $
(respectively $<0$).

We recall from 
H{\"o}rmander \cite[Definition 20.1.1]{Hor:ALPDOIII} (see also \cite[Remark 18.2d]{BooWoj:EBP})
that $P(T)$ defines a {\em local elliptic boundary condition} for the
operator $\wt A$ (or, equivalently, $P(T)$ satisfies the {\v
S}apiro-Lopatinski{\v i} condition for $\wt A$), if and only if the
principal symbol $\sigma^0_{P(T)}$ of $P(T)$ maps the positive spectral
subspace of the leading symbol of the tangential operator, that is
the space $\im P_+(b_0(y,\zeta))\oplus \im (J_0^t)\ii P_-(b_0(y,\zeta)^*)$,
isomorphically onto the fibre $F_y$ for each point $y\in \Sigma$ and
each cotangent vector $\zeta\in T^*_y(\Sigma)$, $\zeta\ne 0$.
This statement is easily seen to be equivalent to
the original definition, which refers to bounded solutions of an ode on the half line.

We now consider the {\v S}apiro-Lopatinski{\v i} mapping $\gs^0_{P(T)}=(-T\, \Id)$ from
$\im P_+(b_0(y,\zeta))\oplus (J_0^t)\ii \im P_-(b_0(y,\zeta)^*)$ to $F_y$:
\[
%\begin{split}
%\im P_+(b_0(y,\zeta))\oplus (J_0^t)\ii \im P_-(b_0(y,\zeta)^*)&\too  F_y\\
\bigl(e_+,(J_0^t)\ii e_-\bigr)  \mapsto  -T e_+ + (J_0^t)\ii e_-.
%\end{split}
\]
Multiplying by $J_0^t$ we see that
this map is bijective if and only if the map
$E_y=\im P_+(b_0(y,\zeta))\oplus \im P_-(b_0(y,\zeta)^*)\too  E_y,$
\begin{equation}\ifarxiv{\label{eq:shap-lop-map}}{\label{ieq:shap-lop-map}}
%\begin{split}
    (e_+,e_-)  \mapsto  -J_0^tT e_+ +e_-
%  \end{split}
\end{equation}
is bijective.
\details{ommitted explaining paragraph:\par

To explain why
$E_y=\im P_+(b_0(y,\zeta))\oplus \im P_-(b_0(y,\zeta)^*)$ we note
that $\im P_+(b_0(y,\zeta))^\perp=\ker P_+(b_0(y,\zeta))^*=\ker P_+(b_0(y,\zeta)^*)
=\im P_-(b_0(y,\zeta)^*)$, so the sum on the left of \eqref{eq:shap-lop-map}
is indeed an orthogonal decomposition.}
%\begin{sloppypar}
Since the dimensions on the left and on the right coincide it
suffices to show injectivity: so let
$-J_0^tT e_++e_-=0,$
$e_+\in \im  P_+(b_0(y,\zeta)),$ $e_-\in \im P_-(b_0(y,\zeta)^*)
=\im P_+(b_0(y,\zeta))^\perp$.
Taking scalar product with $e_+$ we find $0=-\scalar{J_0^tT e_+}{e_+}$.
This implies, since by assumption $J_0^tT>0$, that
$e_+=0$. But then $e_-=0$ as well.

Thus it is proved that $P(T)$ satisfies the
{\v S}apiro-Lopatinski{\v i} condition. This implies that $\wt A_T$ is
a Fredholm operator with compact resolvent. In particular $\wt A_T$
has closed range and hence $\operatorname{coker} \tilde A_T
\cong \ker \tilde A_T^*$.

If $f_+\in Z_0(A), f_-\in Z_0(A^t)$ then
$(f_+, f_-)=f \in \ker \tilde{A}_{T}$ since
${f_-}|_{\pl M}=0=T {f_+}|_{\pl M}$.
Conversely, let
$(f_+, f_-)=f \in \ker \tilde{A}_{T}$. Then certainly
$f_+\in \ker A, f_-\in \ker A^t$ and  ${f_-}|_{\pl M} = T {f_+}|_{\pl M}$.
Since $J_0^tT$ is nonnegative and invertible, the
operator $W:= (J_0^tT)^{\half}$ exists and is invertible.
Now Green's formula \eqref{e:green} yields
\begin{equation}\label{eq3.27}
\begin{split}
\|W {f_+}|_{\pl M}\|^2   &=   \scalar{ {f_+}|_{\pl M}}{J_0^tT {f_+}|_{\pl M}}= \scalar{J_0 {f_+}|_{\pl M}}{ {f_-}|_{\pl M}}\\
                &= -\scalar{Af_+}{f_{-}}+\scalar{f_+}{A^tf_{-}} = 0,
\end{split}
\end{equation}
and since $W$ is invertible we find ${f_+}|_{\pl M} =0$ and
${f_-}|_{\pl M} = T{f_+}|_{\pl M} =0$; hence $f\in Z_0(A)\oplus Z_0(A^t)$.
%\end{sloppypar}
%\ifarxiv{\end{proof}}{
\end{pf*}
%}

\section{Applications}\label{s:appl}

%\noi 3.1 \textit{The \Calderon\ projection}.
\paragraph*{\textup{3.1} \textit{The \Calderon\ projection}.}
From Theorem
\ref{t:invdoub} the Calder{\'o}n projection may be constructed in the
usual way. Let $r_\pm({f_+},{f_-}):=f_\pm$ and
$\varrho_\pm({f_+},{f_-}):= {f_\pm}|_{\partial M}$ and denote by
$\varrho^*$ the $L^2$--dual of $\varrho_+$. It is well known that
$\varrho_\pm$ maps the Sobolev space $L^2_s(M,\ldots)$ continuously
into $L^2_{s-1/2}(\partial M,\ldots)$ for $s>1/2$ and consequently
$\varrho^*$ maps $L^2_s(\partial M,\ldots)$ continuously into
$L^2_{s-1/2}(M,\ldots)$ for $s<0$. These constraints on $s$ cause
some technical difficulties.

For the domain of $A_{\max}$ the mapping properties of $\varrho$ can be
slightly improved. Namely, for $s\geq 0$ the trace map extends by continuity to a
bounded linear map
\begin{multline}\ifarxiv{\label{eq3.11}}{\label{ieq3.11}}
\cD(A_{\max,s})=\bigsetdef{u\in L^2_s(M,E)}{Au\in L^2_s(M,F)}\\ \longrightarrow L_{s-\half}^2({\pl M},E|_{\pl M}), \quad s\geq 0,
\end{multline}
here the domain $\cD(A_{\max,s})$ is equipped with the graph norm of $A$ in $L^2_s(M,E)$.
%that is, there is a constant $C_s$, such that for $f\in L_s^2(M,E)$ with $Af\in L^2_s(M,E)$
%
%\begin{equation}\label{eq3.12}
%\|\varrho f\|_{s-\half} \leq C_s(\|f\|_s+\|Af\|_s)
%\end{equation}
%Here $\|f\|_s$ denotes the Sobolev norm of order s.

\begin{dfn}\label{d:cal}
Let $\wt A_T^{-1}$ denote the pseudoinverse of the operator $\wt
A_T$. That is on $\im \wt A_T$ it is the inverse of $\wt A_T|_{(\ker\wt A_T)^\perp}$,
on $(\im \wt A_T)^\perp$ it is defined to be zero. Put
\[K_\pm:=\pm r_\pm \tilde A_T^{-1}\varrho^*J_0,\quad
C_+:=\varrho_+K_+,\quad C_-:=T^{-1}\varrho_-K_-\,.
\]
\end{dfn}

\begin{theorem}\label{t:calproj} Under the assumptions of Theorem
\textup{\ref{t:invdoub}} and under the additional technical assumption
that the commutator $[B_0^t,J_0^*T]$ is of order $0$ we have:
%\begin{enumerate}

\textup{1.} For $s\ge -1/2$ the operator
$K_+$ maps $L_s^2(\partial M,E|_{\partial M})$ continuously into
$L^2_{s+1/2}(M,E)\cap \ker A$.

\textup{2.} $C_\pm$ are complementary idempotents with
%% \begin{equation*}%% \begin{split}
$\operatorname{im}{C_+}=N_+^0$ and
$\operatorname{im}{C_-}=T^{-1}N_-^0$, $N_-^0$ being the Cauchy data space of $A^t$.
%% \end{split}
%% \end{equation*}
%\end{enumerate}
If $T=(J_0^*)^{-1}$ then $C_\pm$ are orthogonal projections.
\end{theorem}
\begin{remark}\label{r:comment-commutator}
$[J_0^tT,B_0^t]=0$ for the choice $T=(J_0^t)\ii$.

If $A=A^t$ and $B_0-B_0^t$ is of order $0$ then
$[J_0^tT,B_0^t]$ is of order $0$ for
$T\in \bigl\{(J_0^*)\ii,{J_0},{J_0}|J_0|^{-1}\bigr\}$.
\end{remark}
\begin{pf*}{Sketch of proof;
for details cf. \mbox{\cite[pp. 12--20, 23--30]{BooLesZhu:CPD}}}
The method is more interesting than this result, which looks pretty
similar to what one gets from geometric invertible double
constructions respectively invertible non--\-ca\-noni\-cal closed
extensions.

There is one tricky point, however, in 1.
%Namely the constraint
%on $s$ dictated by the Trace Theorem for Sobolev spaces implies
%that $\varrho^*$ maps $L²_s(\pl M,F)$ continuously into
%$L²_{s-1/2}(M,F)$ only if $s<0$. From this the claim in 1. can
%easily be deduced for $-1/2\le s<0$.
Namely the constraint on $s$ dictated by the Trace Theorem for Sobolev spaces,
as explained above. From this the claim in 1. can
easily be deduced only for $-1/2\le s<0$.

To extend it to $s\ge 0$ (including the interesting case $s=0$)
one could invoke the general theory of elliptic boundary value
problems (e.g. Grubb \cite{Gru:FCP}).

We prefer a more elementary approach which is also better suited
to deal with parameter--dependence, see Subsection 3.3 below.

A crucial observation is that the tangential
operator $B_0$ is not an arbitrary elliptic operator. Rather the
ellipticity of $A$ implies that $B_0-it, t\in \mathbb{R},$ is
elliptic in the parametric sense. This is much stronger than
ellipticity (cf. Shubin \cite{Shu:POS} for definition and basic properties).

This observation allows us to introduce the operators
\begin{align}
Q_+(x)&:=\frac 1{2\pi i}\int_{\gG_+} e^{-x \gl}(\gl-B)\ii\, d\gl\,, \qquad x> 0, \ifarxiv{\label{eq:qplus}}{\label{ieq:qplus}}     \\
Q_-(x)&:=\frac 1{2\pi i}\int_{\gG_-} e^{-x \gl}(\gl-B)\ii\, d\gl, \qquad x < 0. \ifarxiv{\label{eq:qminus}}{\label{ieq:qminus}}
\end{align}
Here $\Gamma_+$ is a
contour which encircles the eigenvalues of $B_0$ in the right half
plane and such that $\Re z_n\to \infty$ if $z_n$ is on $\Gamma_+$ with
$|z_n|\to \infty$;  the contour $\Gamma_-$ in the left half plane encircles the complementary
set of eigenvalues, see Figure 1. % \ref funktioniert nicht, liefert Figure 3 (!) 

\begin{figure}
\ifarxiv{% arxiv does not like the nice graphics
\centerline{\includegraphics[height=6cm]{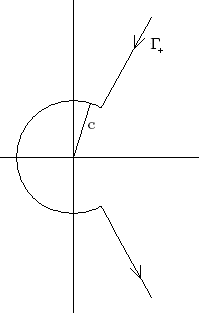}\hspace*{1cm}
\includegraphics[scale=0.6]{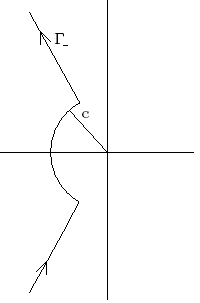}}
}{%
\centerline{\input{fig2.pdf_t}}
}% ifarxiv
\caption{\label{f:gamma_pm} The contours $\gG_{\pm}$ in the plane defining the sectorial
projections $Q_{\pm}$}
\end{figure}
It can be shown (cf. \cite[Sec. 3.2]{BooLesZhu:CPD}) that $P_\pm(B)=Q_\pm(0)$ is an a priory
unbounded idempotent. It follows, however, from the work of Seeley \cite{See:CPE}, Burak \cite{Bur:SPE}, Wojciechowski
\cite{Woj:SFG}, see also Ponge \cite{Pon:SAZ}, that $P_\pm(B)$ is a pseudodifferential operator of order $0$.
Intuitively, the positive/negative sectorial projections $P_\pm$
should map onto the subspaces spanned by the generalized eigenvectors corresponding
to the eigenvalues encircled by the contours $\Gamma_\pm$. See, however,
Remark \plref{warning-root-vectors} below.

The following approximation to $K_+$ constructed from $Q_+$ allows to control
the error when replacing $C_+$ by the sectorial projection of $B_0$.

By direct calculation one checks (\cite[Prop. 3.16]{BooLesZhu:CPD}, cf. also 
Himpel, Kirk and Lesch \cite[Prop. 3.13]{HimKirLes:CPH} where the Spectral Theorem
is used in the case of a symmetric tangential operator) that
for $s\in\R$,  a cut--off function $\varphi\in C^\infty_0(\R_+)$ and $m\in\Z_+$ the operator
\begin{equation}\ifarxiv{\label{eq:nine}}{\label{ieq:nine}}
\id_{\R_+}^m \varphi Q_+: \xi\mapsto \bigl( x \mapsto x^m\/ \varphi(x)\/ Q_+(x)\xi\bigr)
\end{equation}
maps $L^2_s({\pl M},E|_{\pl M})$
continuously to $L^2_{\comp}(\R_+,L^2_{s+m+\half}({\pl M},E|_{\pl M}))$.
Furthermore, for $s\ge -\half$ it maps continuously to\ifarxiv{\\}{\relax}
  $L^2_{s+m+\half,\comp}(\R_+\times{\pl M},E|_{{\pl M}})$. Put
\begin{equation}\label{eq4.20}
    (R\xi)(x):=\begin{pmatrix} Q_+(x) \xi \\
                          Q_-(-x)^*\xi
           \end{pmatrix}, \quad R_T\xi:=\begin{pmatrix} \Id& 0\\ 0 &
    -T \end{pmatrix} R\xi.
\end{equation}
One checks the identity
\begin{equation}\ifarxiv{\label{eq4.21}}{\label{ieq4.21}}
     \wt A_{T}\varphi R_T\xi=
     \bigl(\varrho^*J_0(P_++P_-^*)+S(A,T)\bigr)\xi.
\end{equation}
Beware that $\varphi R_T$ does not map into $\dom \wt A_{T}$.
Therefore \eqref{eq4.21} has to be taken with a grain of salt.
It is an identity in the ``Sobolev space''
$H_{-1}((\wt A_{T})^*)$ which is the dual space of $\dom (\wt A_{T})$,
where the latter is equipped with the graph norm. This fact is reflected
in the appearance of $\varrho^*$. If one applies just the differential
expression $\tilde A$ to $\varphi R_T\xi$ in $M\setminus {\pl M}$ then one obtains,
by definition, $S(A,T)\xi$.

\details{Removed paragraph, duplicates stuff from above....
$P_\pm$ denote
the positive/negative sectorial projections of $B_0$ introduced above.
%$P_+=s-\lim_{x\to 0+}Q_+(x)$
%($P_-=s-\lim_{x\to 0-}Q_-(x)$) is the strong limit as $x\to 0$ of $Q_+(x)$.
A priori these are unbounded idempotents. However, since $B_0+it, t\in\R$,
is elliptic in the parametric sense it follows that $P_\pm$ are
pseudodifferential operators of order $0$ (Burak \cite{Bur:SPE}, Ponge \cite{Pon:SAZ}).}

The reader should be warned that from the equality \eqref{eq4.21} one should not
draw false conclusions about the mapping properties of $\varrho^*$.
$\varphi R_T$ maps $L^2$ to $L^2_\half$ but not to
$H_{\half}(\wt A_{T})$. So one cannot conclude (and it is
indeed not true in general) that the right hand side of
\eqref{eq4.21} lies in $H_{-\half}(\wt A_{T})$.
Nevertheless with some care one can derive the following identity
from \eqref{eq4.21}
\begin{equation}\ifarxiv{\label{eq4.22}}{\label{ieq4.22}}
%\begin{split}
     \wt A_T^{-1} \varrho^*=\Bigl((\Id-P_{Z_{0}(A)})\varphi R_T - \wt A_T^{-1} S(A,T)\Bigr)(J_0(P_++P_-^*))^{-1}.
%\end{split}
\end{equation}

This identity yields the approximation $\varphi R_T ((P_++P_-)^*)^{-1}$ to $K_+$ near the boundary mentioned above.
The mapping properties of $K_+$ can now be derived from those of $\varphi R_T$ and those
of $S(A,T)$.

The mapping properties of $\varphi R_T$ were already quoted above.
To calculate $\wt A \varphi R_T$ we proceed by component:
\begin{equation}\ifarxiv{\label{eq4.18}}{\label{ieq4.18}}
   A\varphi Q_+(x) \xi = \Bigl((C_1 x+C_0)\varphi(x) +J_0
   \varphi'(x)\Bigr)Q_+(x) \xi
\end{equation}
and
\begin{equation}\ifarxiv{\label{eq4.19}}{\label{ieq4.19}}
   \begin{split}
   A^t \varphi T Q_-(-x)^*\xi&=
      \Bigl((\wt C_1 xT+\wt
      C_0T+[B_0^t,J_0^tT])\varphi(x)-\ldots\\
      &\quad  -J_0^t\varphi'(x)\Bigr)Q_-(-x)^* \xi.
\end{split}
\end{equation}
An immediate consequence of \eqref{eq4.18}, \eqref{eq4.19} and \eqref{eq:nine}
is that $S(A,T)$ maps $L^2_s({\pl M},E_{{\pl M}})$ continuously to
$L^2_{s+\half,\comp}(M,F\oplus E)$, $s\ge -\half$.

This establishes the regularity part of the claim. That $K_+$ maps indeed into $\ker A$
follows easily from the fact that $\varrho^*\xi$ is supported on ${\pl M}$ for any $\xi\in L^2_s({\pl M},...)$.

Since $K_+$ maps into $\ker A$ we have  $C_+(L^2)\subset  N_+^0, C_-(L^2)  \subset  T^{-1}N_-^0$.
We show

(i) $N_+^0 \cap T^{-1}N_-^0  =  \{0\}$,

(ii) $C_+ \, +\, C_-  =  \Id$.

This easily implies part 2. of the Theorem.

\enum{i} Let $\xi\in N_+^0 \cap T^{-1}N_-^0$. Then there are
$f\in \ker A\cap L^2_{1/2}(M,E),  g\in \ker A\cap L^2_{1/2}(M,F)$ with $\varrho f = \xi =
T^{-1}\varrho g$. Then
\begin{equation}
(f, g)\in\ker\wt{A}_{T} = Z_0(A)\oplus Z_0(A^t), \quad \text{by Theorem \ref{t:invdoub}.}
\end{equation}
Since elements of $Z_0(A)\oplus Z_0(A^t)$ vanish on the boundary we infer $\xi =
0$.

\enum{ii} Let $\xi\in L^2({\pl M},E|_{\pl M})$,
and $f\in \dom(\wt{A}_{T}^*)$.  Then $\varrho_+ f =-T\varrho_- f$
and exploiting the self--adjointness of $J_0^tT$
we obtain
\begin{equation}
\begin{split}
    \scalar{(C_+&+C_-)\xi}{J_0^t\varrho_+f}
            = \scalar{\varrho_+\wt A_T^{-1}\varrho^* J_0\xi -  T^{-1}\varrho_-\wt A_T^{-1}\varrho^* J_0\xi}{J_0^t\varrho_+ f}\\
        &= \scalar{\varrho_+\wt A_T^{-1}\varrho^*J_0\xi}{J_0^t\varrho_+f}
            - \scalar{\varrho_-\wt A_T^{-1}\varrho^*J_0\xi}{J_0J_0^{-1}(T^{-1})^tJ_0^t\varrho_+f}\\
        &= \scalar{(\varrho_+\oplus\varrho_-)\wt A_T^{-1}\varrho^*J_0\xi}%
       {(J_0^t\oplus J_0)(\varrho_+f\oplus\varrho_-f)}\\
        &= \scalar{\wt A_T^{-1}\varrho^*J_0\xi}{\wt{A}^*_{T}f} = \scalar{\varrho^*J_0\xi}{f}=\scalar{\xi}{J_0^t\varrho f}.
\end{split}
\end{equation}
This proves (ii).

For $T = (J_0^t)^{-1}$ one easily shows using Green's formula
that $N_+^0 \perp T^{-1} N_-^0$, proving that $C_\pm$ are orthogonal projections
in this case.

%To prove the pseudodifferential property we recall from
%\cite[Appendix]{See:TPO} Seeley's construction of the \Calderon\
%projection which always yields a pseudodifferential projection
%$\cP_+^0$ onto $N_+^0$. By orthogonalization within the class of
%pseudodifferential projections and without changing the principal
%symbol and its range (like in \cite[Lemma 12.8]{BooWoj:EBP}), we
%obtain an orthogonal pseudodifferential projection onto $N_+^0$
%which must coincide with $C_+$ for $s=0$.
\end{pf*}
\smallskip

%\noi 3.2 \textit{A general Cobordism Theorem}.
\paragraph*{\textup{3.2} \textit{A general Cobordism Theorem}.}
We shall now give
a wide generalization  of the Cobordism Theorem, previously known
only for operators of Dirac type.

\begin{theorem}[The General Cobordism Theorem]\label{t:cobord}
Let $A$ %% :C^\infty(M,E)\to C^\infty(M,F)$
be a first order formally
self--adjoint elliptic differential operator on a smooth compact
manifold $M$ with boundary acting between sections of the vector bundle
$E$. Then we have the following results:

\enum{I} Let $C_{\pm}$ denote the \Calderon\ projections of
Definition \plref{d:cal}, constructed from the invertible double
with ${T}\in \bigl\{(J_0^*)\ii,{J_0},{J_0}|J_0|^{-1}\bigr\}$.
Then the range of $C_+$ is a Lagrangian subspace of the strongly
symplectic Hilbert space $\bigl(L^2(\partial M,E|_{\partial M}),
-{J_0}\bigr)$. Note that $\im C_+$ is independent of ${T}$.
Consequently, there exists a self-adjoint pseudodifferential
Fredholm extension $A_P$.

\enum{II} We have $\sign i P_0 {J_0}|_{W_0} =0$. Here $W_0$ denotes
the (finite--dimen\-sional) sum of the generalized eigenspaces of
${B_0}$ to imaginary eigenvalues and $P_0$ denotes the orthogonal
projection onto $W_0$; in general ${J_0}$ will not map $W_0$ into
itself. If ${B_0}={B_0}^t$\,, then ${J_0}$ anticommutes with ${B_0}$
and we have $\sign i{J_0}|_{\ker {B_0}} =0$ and the tangential
operator ${B_0}$ is odd with respect to the grading given by the
unitary operator $\ga:=i {J_0}(-{J_0}^2)^{-1/2}$ and hence splits
into matrix form ${B_0}=\bigl(\begin{smallmatrix} 0& B^-\\  B^+&
0\end{smallmatrix}\bigr)$ with respect to the $\pm 1$--eigenspaces
of $\ga$. The index of $B^+:\ker (\ga-1)\longrightarrow \ker(\ga+1)$
vanishes.
\end{theorem}

\begin{pf*}{Sketch of proof;
for details cf. \mbox{\cite[pp. 30--41]{BooLesZhu:CPD}}}
The first claim,
expressed in the language of symplectic functional analysis, follows
immediately from our construction of the \Calderon\ projection:
Clearly $(u,v)\mapsto \lla -{J_0}u,v\rra$ is a (strong)
symplectic form for the Hilbert space $L^2(\partial M,E|_{\pl M})$.
Then the range  $\im(C_+)=N_+^0$ is an isotropic
subspace because of Green's formula \eqref{e:green}. It remains
valid for arbitrary formally self-adjoint elliptic operators of
first order and is here  applied to $\ker A$. Here we use the symplectic $T:=
{J_0}|J_0|^{-1}$ to construct $C_\pm$\,.  That choice of $T$ yields a
self-adjoint double $\tilde A_T$ as observed above.
Then also $\im(C_-)=T\ii(N_+^0)$ is an isotropic
subspace. By Theorem \ref{t:calproj}.2, we have $C_+ + C_- = \Id$, so $N_+^0$
and $T\ii(N_+^0)$ make a pair of transversal isotropic subspaces
of $L^2(\partial M,E|_{\pl M})$. Then (I) follows from
the simple algebraic observation, that transversal isotropic subspaces in a
symplectic vector space must be Lagrangian (taken from Boo{\ss} and Zhu \cite[Lemma 1.2]{BooZhu:WSF}).

To derive the second claim we split $L^2(\partial M,E|_{\pl M})$ into
the spectral subspaces $W_<\,,W_0\,, W_>$ of $B_0$ corresponding to eigenvalues
with negative, respectively zero, respectively positive real parts.
The projections onto $W_<, W_>$ are in fact versions of the positive/negative
sectorial projections which are defined via contour integrals, cf. the proof of Theorem
\ref{t:calproj}.

Then $(\im C_+,W_<)$ is a Fredholm pair by
finite--dimensional perturbation. Hence $\im C_+\oplus W_<$ is a
closed subspace.

We notice that the annihilator $W_<^o=W_0\oplus W_<$ is a co-isotropic subspace of
$L^2(\partial M,E|_{\pl M})$. We
apply symplectic reduction to the Lagrangian subspace $\im C_+$ and obtain that
$
\bigl((\im C_+) + W_<^{oo})\cap W_<^o\bigr)/W_<^{oo}
$
is a Lagrangian subspace of
$
W_<^o/W_<^{oo} = (W_<+W_0)/W_< \simeq W_0.
$
So, the finite--dimensional symplectic Hilbert space
$\bigl(W_0,\scalar{i P_0{J_0}\cdot}{\cdot}\bigr)$
has a Lagrangian subspace. Therefore $\sign i P_0 {J_0}|_{W_0}=0$.
The remainder of (II) follows like in \cite[Theorem 21.5]{BooWoj:EBP}.
\end{pf*}

%% \smallskip

%\noi 3.3 \textit{Parameter dependence}.
\paragraph*{\textup{3.3} \textit{Parameter dependence}. }
To study continuous
variations, we equip the space $\cE(M;E,F)$ consisting of pairs
$(A,T)$ with $J_0^*T>0$ and $[B_0^t,J_0^*T]$ of order 0 with the
metric $d_0((A,T),(A',T'))$ := $N_0(A-A',T-T')$ and the strong
metric $d_{\str}((A,T),(A',T')) := N_0(A-A',T-T') + N_1(A-A',T-T')$,
where
\begin{align*}
N_0(A,T)&:= \|A\|_{1,0}+\|A^t\|_{1,0}+\|T\|_{\half,\half} \quad \tand\\
N_1(A,T)&:=\|B_0\|_{1,0}+\|B_0^t\|_{1,0}+\|[B_0^t,J_0^*T]\|_{0}
                       +\|T\|_{0} \ifarxiv{\label{eq:strong-norm}}{\label{ieq:strong-norm}}\\
                   &\quad +\|J_0\|_{0}+ \|C_1\|_{1,0}+\|C_0\|_0+
                   \|\wt C_1\|_{1,0}+\|\wt C_0\|_0 \,.
\end{align*}
Here $\|\cdot\|_{s,t}$ denotes the norm for bounded operators
from the Sobolev space $L^2_s(...)$ into the Sobolev space $L^2_t(...)$.

We denote by $\cE_{\operatorname{UCP}}(M;E,F)$ the subspace
consisting of pairs $(A,T)$ where $A$ and $A^t$ satisfy weak inner
\UCP\ and by $\cE_{\operatorname{UCP}}^{\sa}(M;E,F)$ the subspace
where $A$ has tangential operator with self--adjoint principal
symbol. We denote by
$\operatorname{Ell}^{\sa}_{\operatorname{UCP}}(M;E)$ the component
of $\cE^{\sa}_{\operatorname{UCP}}$ of formally self--adjoint
operators, equipped with the strong metric.

\begin{theorem}\label{t:param}
\textup{(I)} The map \[(\cE_{\operatorname{UCP}},d_0)\longrightarrow
\cB(L^2(M,F\oplus E),L^2_1(M,E\oplus F)),\quad (A,T)\mapsto \tilde
A_T^{-1}\]
is continuous.

\noi \textup{(II)} For $s\in [-\half,\half]$ the map
\[
(\cE^{\sa}_{\operatorname{UCP}}(M;E),\dstr)\longrightarrow
\cB(L^2_s(\partial M,E|_{\partial M})),\quad  (A,T) \mapsto C_+(A,T)
\]
is continuous.

\noi \textup{(III)} The map $
\operatorname{Ell}^{\sa}_{\operatorname{UCP}}(M;E) \longrightarrow
\cB(L^2_1(M,E),L^2(M,E))$, $A\mapsto A_{C_+}$ is continuous. Here
$C_+$ denotes the version of the \Calderon\ projection constructed
from $T:=(J_0^*)\ii$.
\end{theorem}

\begin{pf*}{Sketch of proof;
for details cf. \mbox{\cite[pp. 42--51]{BooLesZhu:CPD}}}

\noi To (I):
The difficulty we are facing here is that
$\dom \bigl( \wt A_T\bigr)$  varies with $T$. So
we fix a $T_0$ close to $T$ and make the following
factorization
\begin{equation}\ifarxiv{\label{e:factorization}}{\label{ie:factorization}}
(A,T)\mapsto \tilde A_{T}\circ \Phi_{T_0,T}
\mapsto \bigl(\tilde A_T\circ
\Phi_{T_0,T}|_{\dom(\tilde A_{T_0})}\bigr)^{-1}
\mapsto \tilde A_T\ii\,.
\end{equation}
Here we set
\[
    \Phi_{T_0,T}\binom{f_+}{f_-}:=\binom{f_+}{f_-+e(T-T_0)\varrho
    f_+},
\]
where $e$ denotes a linear right-inverse to $\varrho$ (as, e.g., in
\cite[Definition 11.7e]{BooWoj:EBP}).
The operator $\Phi_{T_0,T}$ is bounded invertible on $L^2_1(M,E\oplus F)$
and maps $\dom \tilde A_{T_0}$ bijectively onto $\dom \tilde A_{T}$.
Note that $\tilde A_T\ii = \Phi_{T_0,T}\circ
\bigl(\tilde A_T\circ \Phi_{T_0,T}|_{\dom(\tilde A_{T_0})}\bigr)^{-1}$\,.
It is straightforward to check the continuity of the first and last arrow
in \eqref{e:factorization}, while the middle arrow is just the general continuous
inversion for bounded operators in Banach space.

\noi To (II): From \eqref{eq4.21} and \eqref{eq4.22} we obtain the correction formula
\begin{equation}\ifarxiv{\label{eq4.24}}{\label{ieq4.24}}
      C_+= \bigl(P_+ - \varrho_+ \tilde A_T^{-1} S(A,T)\bigr)(P_++P_-^*)\ii\,,
\end{equation}
relating $C_+(A,T)$ and the sectorial projection $P_+:=P_+(B_0)$ of the
tangential operator $B(0)$.
Note that \eqref{eq4.24} is valid for arbitrary $A$. Continuous
variation of $\tilde A_T^{-1}$ and $S(A,T)$ can be obtained without restriction regarding $B_0$
in the same way as in the proof of (I). Continuous variation of $P_+$ can also be obtained
for general $A$ and general $B_0$ if the variation is of order $<1$ \cite[Prop. 7.13]{BooLesZhu:CPD}.
However, admitting
variation of $A$ and $B_0$ of order 1, we obtain continuous variation of $P_+$ only for
formally self-adjoint $B_0$ via a Riesz-map argument,
see the discussion around Prop. \plref{p:dependence-Pplus-self-adjoint} below.

\noi To (III): Note that $A_{C_+}$ is self--adjoint by Theorem \ref{t:cobord}.II.
It is indeed a self--adjoint realization of
a well--posed boundary value problem.
Applying the preceding (II) and a generalization of the preceding (I)
(as indicated below in the Note to (I)) yields the claim.
\end{pf*}

\begin{remark} (I) is much stronger than just graph continuity. In
the same way, we obtain that the map
\[
(A,P)\mapsto (A_P+i)\ii\in
\cB(L^2(M;E),L^2_1(M;E))
\]
is continuous with respect to the $d_0$
metric on the space of pairs $(A,P)$ where $P$ is a
pseudodifferential orthogonal projection which is well--posed with
respect to $A$. In particular $(A,P)\mapsto (A_P+i)\ii$ is graph
continuous.

\noi The continuous dependence of the Calder{\'o}n projection on the input
data in (II) has consequences for the so--called Spectral Flow Theorem, cf.
\cite{BooZhu:WSF} and Section 1.4.B above.

\noi In (III) we obtain a more precise version of
\cite{BooLesPhi:UFOSF}, Theorem 3.9 (c). Note that our present
version applies to a much wider class of operators than loc. cit.

For deformations of $A$ by a $0$th order operator the 
continuous variation of the \Calderon\ projector was proved
in a purely functional analytic context by Boo{\ss}--Bavnbek and Furutani \cite{BooFur:MIF}.
\end{remark}

\section{Closing, respectively confining, the gaps}\label{s:gaps}

\paragraph*{\textup{4.1} \textit{Full elucidation of parameter dependence.}}
Dirac type operators and their
\Calderon\ projection vary continuously under variation of the underlying connection.
Note that such variations are only perturbations by bounded operators. However,
it is a disturbing gap in our perception of
fundamental concepts of quantum field theory that we can not exclude the possibility that
more general but still smooth variations of the coefficients can yield jumps of the
\Calderon\ projection when leaving the realm of Dirac type operators, 
even in the presence of a spectral cut.

%% \subsection{General variations}

In particular, for spectral flow formulas (see Subsection 1.4.B),
the continuous dependence of the \Calderon\ projection on the input data is crucial.
In Section 3 we could establish this continuous dependence under the assumption
that the tangential operator has self-adjoint leading symbol.
We would like to get rid of this technical assumption.
For this it is necessary to study the positive sectorial projection of a (non-self-adjoint) elliptic
differential operator and its dependence on the input data.

\details{Removed paragraphs (ML 15.09.2008):\par\medskip

First we address open questions regarding the relation between the metrics $d_0$ and $d_{\str}$ and 
the norms $N_0$ and $N_1$ introduced above in Subsection 3.3. Basically, that is the question 
of the continuous or not--continuous variation of the tangential operator and the principal symbol
in normal direction along the boundary when the operator as a whole varies continuously.

Then we recall the power of symbolic calculus and arguing modulo operators 
of lower order in chains of pseudodifferential operators, as in Seeley \cite{See:CPE}.

Finally, we explain its limitation for studying continuous parameter dependence (like in
email correspondence with Chaofeng Zhu).

}% end of removed paragraphs

In \cite[Theorem 7.2]{BooLesZhu:CPD} we in fact proved the following variant of Theorem \plref{t:param}
(II) for $(\cE_{\operatorname{UCP}}(M;E),\dstr)$ instead of $(\cE^{\sa}_{\operatorname{UCP}}(M;E),\dstr)$:
\begin{sloppypar}
If $(A(z),T(z))$ is a continuous family in $(\cE_{\operatorname{UCP}}(M;E),\dstr)$ such that
the corresponding family of positive sectorial projections of the tangential operator
varies continuously in $\cB(L^2_s(\partial M,E|_{\partial M}))$ then 
$C_+(A(z),T(z))\in \cB(L^2_s(\partial M,E|_{\partial M}))$ depends continuously on $z$, too.
\end{sloppypar}

So the continuous dependence of the \Calderon\ projection in general 
hinges on the continuity of the positive sectorial projection of the tangential operator. 

We can now phrase the problem completely in terms of a single elliptic operator on the closed
manifold $\partial M$: fix contours $\Gamma_+, \Gamma_-$ as in Figure 1 and denote by % \ref funktioniert nicht, liefert Figure 3 (!) 
$\Ell^1_\Gamma(\partial M,E|_{\partial M})$
the space those differential operators $B$ on $\partial M$ acting on sections
of the vector bundle $E|_{\partial M}$ such that 
firstly $B$ is elliptic and all eigenvalues of the leading symbol $\sigma^1_B(p,\xi)$
are encircled by $\Gamma_+ \cup\Gamma_-$ and secondly 
no eigenvalue of $B$ lies on the curves $\Gamma_\pm$.

The first condition means in other words that for $\xi\in T_p^*M\setminus\{0\}$ there are no eigenvalues of the leading symbol of $B$ inside the sectors
around the imaginary axis with boundary defined by the legs of $\Gamma_\pm$. 

Note that a self--adjoint operator is in $\Ell^1_\Gamma(\partial M,E|_{\partial M})$ if and only if
$\pm c\not\in\spec B$. 

For $B\in\Ell_\Gamma^1(\partial M,E|_{\partial M})$ the operator families $Q_\pm$, in particular the positive
sectorial projections $P_\pm$, as defined in the proof of Theorem \ref{t:calproj}, are available.

Now we can formulate the main problem related to the continuous variation of the \Calderon\ projection:

\begin{problem} \label{p:variation-positive-projection}
For which (reasonable) topology on $\Ell_\Gamma^1(\partial M,E|_{\partial M})$ is the map
\begin{equation}
    \Ell_\Gamma^1(\partial M,E|_{\partial M}) \longrightarrow \cB(L^2(\partial M,E|_{\partial M}), 
\quad B\mapsto P_\pm(B)
\end{equation}
continuous?
\end{problem}

For our original operator $A$ with tangential operator $B_0$ (cf. \eqref{eq:operator-collar})
we only have a reasonable chance for $P_\pm(B_0)$ to depend continuously on $A$ if $B_0$ 
depends continuously on $A$, too. So the next problem is

\begin{problem} Suppose the topology $\alpha$ of Problem \plref{p:variation-positive-projection}
is found. For which topology $\beta$ on $\cE_{\operatorname{UCP}}(M;E)$ is then
the map 
\[
(\cE_{\operatorname{UCP}}(M;E),\beta)\longrightarrow (\Ell^1_\Gamma(M;E),\alpha),\quad 
    A\mapsto B_0(A)\]
continuous?
\end{problem}

The construction of the \Calderon\ projection dictates that the topology $\beta$ must be finer
than the one induced by the strong metric to ensure then the continuity of the
map 
\[(\cE^{\sa}_{\operatorname{UCP}}(M;E),\beta))\ni (A,T)\mapsto C_+(A,T)\in 
\cB(L^2_s(\partial M,E|_{\partial M})).\] \details{Ich kann nicht sehen, we wie argumentieren
soll, dass die $d_0$ Metrik schon genuegen soll. Any suggestions??}

The power of the Spectral Theorem is amazing. Namely, if we restrict to self--adjoint tangential
operators $B_0$ then the Spectral
Theorem allows to prove (cf. \cite[Prop. 7.15]{BooLesZhu:CPD}.

\begin{prop}\label{p:dependence-Pplus-self-adjoint}
Let $\Ell_c^{1,\sa}(\partial M,E|_{\partial M})$ denote the space of \emph{self--adjoint}
first order elliptic differential operators B with $\pm c\not\in\spec B$.
Then for $|s|\le \half$ the map
\begin{equation}
\begin{split}      
\Bigl(\Ell_c^1,\|\cdot\|_{1,0}\Bigr)&\longrightarrow \cB(L^2_s(\partial M,E|_{\partial M}))\\
              B&\mapsto 1_{[c,\infty)}(B)
            \end{split}
\end{equation}
is continuous.
\end{prop}

As a consequence for $\cE^{\sa}_{\operatorname{UCP}}(M;E)$ we can choose $\beta$ to be the
topology induced by the strong metric $d_{\str}$ and $\alpha$ to be the topology induced
by the norm $\|\cdot\|_{1,0}$. 

So besides the obvious fact that eigenvalues crossing $\pm c$ (resp. the contours of integration in
the non--self--adjoint case) lead to jumps of $P_+(B)$, for self--adjoint $B$ 
the map $B\mapsto P_+(B)$ is even continuous when $\Ell^1_c(\partial M,E|_{\partial M})$
is equipped with the relatively weak norm topology of bounded maps from $L^2_1(\partial M,E|_{\partial M})$ to 
$L^2(\partial M,E|_{\partial M})$.

For general non--self--adjoint $B$ we are far from being able to prove such a result.

Let us focus on $P_+(B)$ in Problem \plref{p:variation-positive-projection}. 
The problem is subtle since the definition of $P_+(B)$
is a bit tricky. A priori it is defined via the ill--defined integral
\begin{equation}\ifarxiv{\label{eq:illdefined}}{\label{ieq:illdefined}}
     P_+(B)=\frac 1{2\pi i}\int_{\gG_+} (\gl-B)\ii\, d\gl
\end{equation}
which does not converge and hence does not easily allow norm estimates.
It does not help much that for $\xi\in\dom(B)$ 
the operator $P_+(B)$ is given by the well--defined integral
\begin{equation}\ifarxiv{\label{eq:welldefined}}{\label{ieq:welldefined}}
     P_+(B)\xi=\xi+\frac 1{2\pi i}\int_{\gG_+} \gl\ii(\gl-B)\ii\, d\gl (B\xi).
\end{equation}
The latter representation is not better suited for
proving norm estimates of the form $\| P_+(B)-P_+(\tilde B)\|$ 
in $L^2$ because the $L^2$--norm of \eqref{eq:welldefined} is a priori unbounded for arbitrary $\xi\in L^2$.

So let us equip $\Ell^1_\Gamma(\partial M,E|_{\partial M})$ with the natural Fr{\'e}chet topology
on (pseudo)differential operators and try to exploit the power of the symbolic calculus:

Let $\Omega$ be an open set in the plane containing a conic neighborhood of the legs of the
contours $\Gamma_\pm$ and containing the closed disc of radius $c$. We assume that $c$ and $\Omega$
are deliberately chosen such that for any cotangent vector $\xi\in T_p^*M$ of length $\ge 1$ the spectrum
of the leading symbol $b(x,\xi):=\sigma^1_B(p,\xi)$ does not meet $\Omega$. 
The resolvent $(B-\gl)\ii\in\CL^{-1}(M;E)$ is now a classical pseudodifferential operator of order $-1$
in the parametric calculus with parameter $\gl\in\Omega$. Choose a cut--off function $\psi\in C^\infty(\R)$
with $\psi(u)=0$ for $|u|\le 1$ and $\psi(u)=1$ for $|u|\ge 2$. Then
$r_{-1}(x,\xi;\gl)=\psi(\sqrt{|\xi|^2+|\gl|^2})(b_0(x,\xi)-\gl)\ii$ is a symbol of order $-1$ in the parameter dependent
calculus and we may consider $R(\gl):=\Op(r_{-1})$. $R(\gl)$ approximates $(B-\gl)\ii$ up to operators of
order $-2$ and therefore (cf. \cite[Theorem 9.3]{Shu:POS})
\begin{equation}
           \|(B-\gl)\ii - R(\gl)\|_{L^2}=O(\gl^{-2}), \quad \gl\to\infty \text{ in } \Omega.
\end{equation}
Then certainly 
\begin{equation}\ifarxiv{\label{eq:ill-modified}}{\label{ieq:ill-modified}}
     \frac 1{2\pi i}\int_{\gG_+} \bigl((\gl-B)\ii\,-R(\gl)\bigr)\, d\gl
\end{equation}
is well-defined. The integral $\frac 1{2\pi i}\int_{\gG_+} R(\gl) d\gl$ can be made sense
of at the symbolic level, where the contour $\Gamma_+$ can be replaced by a closed
contour encircling the eigenvalues of the leading symbol. This construction sketches
why $P_+(B)$ is indeed a pseudodifferential operator of order $0$.

Things become more involved if we assume that everything depends on an additional parameter, say
$s$. Sufficient for the continuity of $s\mapsto P_+(B_s)$ would be to establish the estimates
\begin{equation}\ifarxiv{\label{eq:cont-one}}{\label{ieq:cont-one}}
\| (B_s-\gl)\ii -R_s(\gl)- (B_{s'}-\gl)\ii +R_{s'}(\gl)\|_{0,0}\le \go(|s-s'|) |\gl|^{-2}
\end{equation}
and
\begin{equation}\ifarxiv{\label{eq:cont-two}}{\label{ieq:cont-two}}
    \Bigl\| \int_{\Gamma_\pm}  R_s(\gl)-R_{s'}(\gl) d\gl\Bigr\| \le \go(|s-s'|)
\end{equation} 
for some function $\go:[0,\infty)\to\R$ with $\lim\limits_{u\to 0+} \go(u)=0$.
 
The problem with \eqref{eq:cont-one} is that $R_s(\gl)$ approximates only asymptotically
in $\gl$. The difference $(B_s-\gl)\ii-R_s(\gl)$ is not necessarily small in norm. 
It would not help here if one would dig deeper into the symbol expansion of $(B_s-\gl)\ii$;
this would only improve the order of approximation in $\lambda$.

Secondly the problem with \eqref{eq:cont-two} is that the left hand side is ill--defined because
$R_s(\gl)-R_{s'}(\gl)=O(\gl^{-1})$, which cannot be improved. One could try to proceed as sketched above: at the symbolic
level one can probably replace the unbounded contour by a closed contour in the plane
encircling the eigenvalues of the leading symbol. 

Nevertheless, the details remain cumbersome. We hope to come back to this problem in a future publication.
 
\paragraph*{\textup{4.2} \textit{The regime of validity of weak inner \UCP}.}
In this Note we have removed any assumption about weak inner \UCP\ from our  canonical construction of the \Calderon\ projection (Definition \ref{d:cal}) and from the General Cobordism Theorem (Theorem \ref{t:cobord}). However, parts of our investigation of the parameter dependence of Poisson operator, \Calderon\ projection, and the continuity of families of ``well--posed"
self--adjoint extensions are only valid under the assumption of weak inner \UCP. It seems to us, in particular, that we need that assumption for establishing the continuity of the changes of the \Calderon\ projection under continuous change of the coefficients of the underlying differential operator.

Then, what is the status of weak inner \UCP\ for elliptic differential operators of first order? 
To further clarify the regime of validity
(or non-validity) of weak \UCP\  to the boundary we have two marks. 

\noi  (I) On the positive side, we have our Theorem \ref{t:ucp-symmetric}. 
However, the partial integrations in the proof of inequality \eqref{e:carleman} depend on the symmetry of the tangential operator (or its elliptic symmetrization $\frac 12 ({B_x}+{B_x}^t)$). So, by now it is hard to see how to get rid of the assumption of symmetric tangential principal symbol. However, see also Proposition \ref{p:ucp-constant} and Conjecture \ref{con:1} below.

\noi (II) On the negative side, we have the classic results of Pli{\v s} \cite[Theorem 2]{Pli:SLE},
where a rather intricate smooth perturbation is given of the bi--harmonic equation which has a smooth non--trivial
solution $u$ on $\mathbb{R}^3$ with $\supp u\< B^3$\,. 
In the same paper, but in different context, Pli{\v s} makes the laconic statement (Remark 4), 
that a certain elliptic equation of fourth order ``is equivalent to the system of four complex 
or eight real equations of the first order". Of course, one can always make the usual textbook substitution of
higher derivatives by new variables, familiar from the treatment of ordinary differential equations of higher order.
However, for dimensions greater than 1, one would loose ellipticity by that way. Alternatively, one could
make a factorization $\gD=\mathbb{D}^2$ of the Laplacian by a suitable restriction of the euclidean 
Dirac operator $\mathbb{D}$, like, e.g.,
in B{\" a}r \cite[Example, Equation (1)]{Bae:ZSS} and, hopefully, extend the factorization 
to the perturbation in a suitable way. Nobody has done that, yet. 

We conclude: Inspired by previous work of S. Alinhac, 
there is a related counter--example for strong \UCP\ also for first order elliptic systems in \cite[l.c.]{Bae:ZSS}, 
while we consider it still an open problem under what conditions weak inner \UCP\
is valid or not for a linear elliptic differential operator of first order with smooth coefficients.

Nevertheless, something more can be proved and much more can be conjectured.
Below, we shall explain a side result of our construction of sectorial projections (see above
Subsection 3.1), namely how the well--known uniqueness of solutions of initial problems for 
systems of ordinary differential equations can be transferred to the case of constant coefficients 
in normal direction in a cylindrical collar of a manifold with boundary, thus preserving weak \UCP. 
Moreover, we discuss the stability of weak inner \UCP\ under ``small" perturbations and non-stability under 
``large" perturbations; and we come up with two conjectures: the first suggesting a criterion for weak \UCP\ 
in terms of the range of the positive sectorial projection, the second affirming Laurent Schwartz's conjecture 
of 1956, that the \UCP\ defects for an elliptic operator of first order and its formal adjoint coincide.

\begin{prop}\label{p:ucp-constant}
Let $\gS$ be a closed manifold and let $A=\dd x+B$ be elliptic on $\R_+\times \gS$, where $B$ is a fixed elliptic (not necessarily formally self--adjoint) operator on $\gS$. Let $u\in\Ci([0,T)\times \gS,E)$ be a section with
\[
Au=0 \quad\tand\quad u|_{{\{0\}\times \gS}}=0.
\]
Then $u=0$.
\end{prop}

\begin{remark}\label{warning-root-vectors}
Let us add a warning: it is tempting to write
\[
A=\bigoplus_{\gl\in\spec B} \dd x+\gL,
\]
where $\gL$ denotes a Jordan matrix with diagonal entries $\gl$ and with $\dd x+\gL$ acting on functions with values in the generalized $\gl$--eigenspace of $B$. Then applying the Picard uniqueness for first order ordinary differential equations one would get the result.

This argument, however, is wrong since it is unknown
whether $B$ has a complete set of root vectors meaning that the sum of the generalized eigenspaces is dense.
For a discussion of this issue and its history see \cite[Sec. 3 and Appendix]{Pon:SAZ}.
There are counterexamples of elliptic differential operators
(Seeley \cite{See:SES} and Agranovich-Markus \cite{AgrMar:SPE}) without a complete set of root vectors. 
In these examples, however, the principal symbol does not admit a spectral cutting.
Our $B$ has the imaginary axis as a spectral cutting such that these counterexamples do not apply. 
This is not enough, however, to apply
the known positive result, see \cite[Appendix]{Pon:SAZ} for details.

Fortunately, we can circumvent that difficulty by applying the sectorial projections studied in Section \ref{s:appl}.
\end{remark}

\begin{pf}
Let
%\begin{align}
%Q_+(z)&:=\frac 1{2\pi i}\int_{\gG_+} e^{-z \gl}(\gl-B)\ii\, d\gl\,, \\
%%% \qquad x> 0,\\
%Q_-(z)&:=\frac 1{2\pi i}\int_{\gG_-} e^{-z \gl}(\gl-B)\ii\, d\gl
%% ,  \qquad x < 0.
%\end{align}
$Q_\pm(z)$ be as in \eqref{eq:qplus}, \eqref{eq:qminus} but with complex argument.
Note that $Q_+(z)$ is analytic for $z$ in a sector containing $(0,\infty)$. Obviously there is a small sector $S_+$ containing $(0,\infty)$ such that for $z\in S_+$ and $\gl\in\gG_+$, $|\gl|\gg 0$
we have $\Re(z\gl)>0$.

Similarly there is a small sector $S_-$ containing $(-\infty,0)$ such that $z\mapsto Q_-(z)$ is analytic for $z\in S_-$\,. Finally recall that $Q_\pm(0)$ are complementary idempotents. We shall not need the full strength of Burak's result. We shall only use the elementary fact that
\[
\lim_{x\to 0_+} Q_+(x)u=Q_+(0)u
\]
exists in $L^2$ for $u\in L^2_s\,, s>0$ and that $Q_\pm(0)$ are (possibly unbounded) idempotents with $Q_+(0)+Q_-(0)=\Id$.

So, let $u\in L^2([0,T)\times \gS,E)$ with
\begin{equation}\label{e:constant}
(\partial_x+B)u(x)=0 \quad\tand\quad u(0)=0.
\end{equation}
Extending $u$ by $0$ to $x<0$ we see that $Au=0$ on the whole interval $(-\infty,T)$; this argument is so easy since
$A$ has constant coefficients and therefore extending $A$ to $x<0$ is trivial. 
By ellipticity of $A$ we then conclude that $u$ is smooth, i.e. $u\in \Ci([0,T]\times\Sigma,E)$.

Consider first for $0< t <T$
\[
g(x):=Q_+(t-x)u(x),\quad 0\le x < t.
\]
From
\begin{align*}
g(0)&=0,\\
g'(x)&= BQ_+(t-x)u(x)-Q_+(t-x)Bu(x)=0
\end{align*}
we infer $g(x)\equiv 0$. Thus
\[
Q_+(0)u(t)=0,\quad 0< t< T.
\]
Next consider for $t<0$
\[
g(x):=Q_-(t-x)u(x),\quad 0\le x < T.
\]
As above we conclude $g(x)\equiv 0$. Consequently,
\[
Q_-(t)u(x)=0
\]
for all $t\le -T$ and all $0\le x <T$.

So, for fixed $0< x< T$ the analytic function
$
z\mapsto Q_-(z)u(x)
$
vanishes for $z\in (-\infty,-T]$ and thus vanishes for all $z\in S_-$\,. In particular $Q_-(0)u(x)=0$.

Summing up we have proved
\[
u(x)=Q_+(0)u(x)+Q_-(0)u(x)=0, \quad 0< x <T.
\]
\end{pf}

\begin{remark}\label{r:ucp-constant}
Note that the case $Q_+(0)u(x)=0$ is very similar to the standard uniqueness proof for 
operator semigroups. However, the case $Q_-(0)u(x) = 0$ is difficult 
because here one needs a kind of uniqueness for a {\em backward heat equation}.
\end{remark}

In view of the prominent role of the sectorial projections in the preceding proposition we risk the following

\begin{conjecture}\label{con:1}
There is a criterion for weak \UCP\ for elliptic differential operators of first order in terms
of the range of the positive sectorial projection $P_+({B_0})$.
%% \marginpar{Relate to \eqref{e:plis-reduction}!?}
\end{conjecture}

The idea of the preceding conjecture is to replace the assumption of Theorem \ref{t:ucp-symmetric}, namely symmetric principal symbol of all tangential operators along all hypersurfaces, by the weaker assumption that the range of all corresponding sectorial projections is dense in $L^2$\,.

One ingredient for the wanted proof may be the {\em \UCP--defect dimension}
\begin{equation}\label{e:ucp-deficiency}
d(x):=\dim\bigsetdef{u}{Au=0\tand u|_{\Sigma(x)}=0}.
\end{equation}
Here $\gS(x)=\{x\}\times \partial M$ denotes the parallel hypersurface in the collar in distance $x$ of the boundary.
Clearly, $d(x)$ is upper semi-continuous, more precisely,
decreasing, left-continuous, and $\Z_+$--valued. The essential result of the preceding proposition is that the UCP--defect dimension
$d(x)$ is constant when the coefficients in normal direction are
constant. It remains to see whether similar constancy results can be obtained for variable coefficients.

%% \medskip

Another noteworthy puzzle, mentioned before, is the question raised by Laurent Schwartz in 1956, see \cite{Sch:EDP}: 
if $A$ is not symmetric, then he asks whether weak \UCP\ for the formal adjoint operator $A^t$ follows from weak 
\UCP\ for $A$.
A related problem (and a decisive one, e.g., for the validity of the Bojarski conjecture, see \cite[Chapter 24]{BooWoj:EBP})
is the vanishing of the {\em inner index}
\[
\ind_0A:=d(0)-d'(0), 
\]
where $d'(x):=\dim\bigsetdef{u}{A^tu=0 \tand u|_{\Sigma(x)}=0}$.
Note that $\ind_0A=\ind A_{C_+}$, i.e., the index of the regular elliptic boundary value problem
with domain defined by the \Calderon\ projection $C_+$. Note also the stability of the \Calderon\
projection established in Section \ref{s:appl}.

It seems to us, however, that this stability does not imply the stability of $\ind_0A$ under smooth shrinking of the underlying manifold. As worked out years ago by the first author, such deformation stability would 
yield the vanishing of the inner index. So, we are not afraid to set forth a second conjecture:

\begin{conjecture}\label{con:2}
The inner index vanishes for elliptic differential operators of first order on compact manifolds with smooth boundary.
\end{conjecture}

We close this Subsection with a few additional comments regarding stability of weak inner \UCP. 
In \cite{BooMarWan:WUCP} two types of very simple examples were given, where perturbation of the standard first order ordinary differential equation $\dot u=0$ by a $0$th order term does not preserve \UCP. 
The first type are sufficiently substantial {\em non--linear} perturbations like adding 
$-2\sqrt{|u|}$ or $-3u^{2/3}$\,. The second type are sufficiently substantial 
{\em global linear} perturbations (say on the interval $[0,2]$)
like adding $-a(x)\int_0^2u(s)a(s)ds$, where $a:[0,2]\to\R$ is a continuous function which vanishes on $[0,1]$ and satisfies $\int_1^2a(s)ds=\sqrt{2}$.

For both types the striking feature is that the \UCP--destructive perturbations are not small. 
More precisely, Boo{\ss} and Zhu \cite[Lemma 3.2]{BooZhu:GSF} prove by semi--continuity of the kernel dimension

%% \newpage

\begin{lemma}\label{l:ucp-stability}
Let $H$ be a Hilbert space. Let $A_s\in \mathcal{C}(H)$, $0\le s\le 1$ be a family of symmetric operators with fixed (minimal) domain $\dom A_s=D_m$ and fixed maximal domain $\dom A_s^*=D_{\max}$. Assume that $\{A_s^*:D_{\max}\to H\}$ is a continuous curve of bounded operators, where the norm on $D_{\max}$ is the graph norm induced by $A_0^*$\,. If $A_0$ satisfies inner weak \UCP\ and there exists a self--adjoint Fredholm extension $A_0^*|_D$ of $A_0$\,, 
then for all $s\ll 1$ the operators $A_s^*$ are surjective and the operators $A_s$ satisfy weak inner \UCP.
\end{lemma}

% etc, etc

% The Appendices part is started with the command \appendix;
% appendix sections are then done as normal sections
% \appendix

% \section{}
% \label{i}

% The Acknowledgements are an un-numbered section
%\section*{Acknowledgements}
% Acknowledgements text here
\begin{acknowledgements}
We would like to thank the editor for thoughtful comments and helpful suggestions
which led to many improvements. The editor clearly went beyond the call of duty,
and we are indebted.
\end{acknowledgements}

\ifarxiv{}{\end{article}}
\end{document}